\documentclass[12pt]{amsart}

\usepackage{graphicx, overpic}
\usepackage[below]{placeins}
\usepackage[colorlinks=true, linkcolor=blue, citecolor=blue]{hyperref}
\usepackage[]{algorithm2e}

\usepackage[T1]{fontenc}

\usepackage{amsmath,amsthm,amscd,amssymb,eucal}

\usepackage{enumerate, amsfonts, latexsym, color, url}
\usepackage{epstopdf}

\usepackage{pinlabel}

\makeatletter
\newsavebox{\@brx}
\newcommand{\llangle}[1][]{\savebox{\@brx}{\(\m@th{#1\langle}\)}%
  \mathopen{\copy\@brx\kern-0.5\wd\@brx\usebox{\@brx}}}
\newcommand{\rrangle}[1][]{\savebox{\@brx}{\(\m@th{#1\rangle}\)}%
  \mathclose{\copy\@brx\kern-0.5\wd\@brx\usebox{\@brx}}}
\makeatother

\begin{document}

\newtheorem{theorem}{Theorem}[section]
\newtheorem{theoremm}[theorem]{`Theorem'}
\newtheorem{lemma}[theorem]{Lemma}
\newtheorem{proposition}[theorem]{Proposition}
\newtheorem{corollary}[theorem]{Corollary}
\newtheorem{conjecture}[theorem]{Conjecture}
\newtheorem{question}[theorem]{Question}
\newtheorem{problem}[theorem]{Problem}
\newtheorem*{claim}{Claim}
\newtheorem*{criterion}{Criterion}
\newtheorem*{universal_circle_theorem}{Universal Circle Theorem~\ref{theorem:universal_circle}}
\newtheorem*{quasigeodesic_theorem}{Quasigeodesic Flow Theorem~\ref{theorem:zippers_from_flows}}
\newtheorem*{quasimorphism_theorem}{Uniform Quasimorphism Theorem~\ref{theorem:zippers_from_quasimorphisms}}
\newtheorem*{order_theorem}{Uniform Order Theorem~\ref{theorem:zippers_from_orders}}

\theoremstyle{definition}
\newtheorem{definition}[theorem]{Definition}
\newtheorem{construction}[theorem]{Construction}
\newtheorem{notation}[theorem]{Notation}
\newtheorem{object}[theorem]{Object}
\newtheorem{operation}[theorem]{Operation}

\theoremstyle{remark}
\newtheorem{remark}[theorem]{Remark}
\newtheorem{example}[theorem]{Example}

\numberwithin{equation}{subsection}

\newcommand\id{\textnormal{id}}

\newcommand\N{\mathbb N}
\newcommand\Z{\mathbb Z}
\newcommand\R{\mathbb R}
\newcommand\C{\mathbb C}
\newcommand\EE{\mathcal E}
\renewcommand\H{\mathbb H}
\newcommand\A{\mathcal A}
\newcommand\F{\mathcal F}
\newcommand\G{\mathcal G}
\newcommand\CC{\mathcal C}
\newcommand\PP{\mathcal P}
\newcommand\QQ{\mathcal Q}
\newcommand\Sp{\textnormal{Sp}}
\newcommand\SL{\textnormal{SL}}
\newcommand\Homeo{\textnormal{Homeo}}
\newcommand\inte{\textnormal{int}}
\newcommand\un{\textnormal{univ}}

\title{Zippers}

\author{Danny Calegari}
\address{University of Chicago \\ Chicago, Ill 60637 USA}
\email{dannyc@uchicago.edu}
\author{Ino Loukidou}
\address{University of Chicago \\ Chicago, Ill 60637 USA}
\email{thelouk@uchicago.edu}
\date{\today}

\begin{abstract}
If $M$ is a hyperbolic 3-manifold fibering over the circle, the fundamental group of $M$ 
acts faithfully by homeomorphisms on a circle (the circle at infinity of the 
universal cover of the fiber), preserving a pair of invariant (stable and unstable) 
laminations. Many different kinds of dynamical structures (e.g.\/ taut foliations, 
quasigeodesic or pseudo-Anosov flows) are known to give rise to universal circles --- 
a circle with a faithful $\pi_1(M)$ action preserving a pair of invariant 
laminations --- and these universal circles play a key role in relating the dynamical structure 
to the geometry of $M$. In this paper we introduce the idea of {\em zippers}, 
which give a new and direct way to construct universal circles, streamlining 
the known constructions in many cases, and giving a host of new constructions in others.
In particular, zippers (and their associated universal circles)
may be constructed directly from {\em uniform quasimorphisms} or from {\em uniform actions}.
\end{abstract}

\maketitle
\setcounter{tocdepth}{1}
\tableofcontents

\begin{figure}[htpb]
\centering
\includegraphics[scale=0.3]{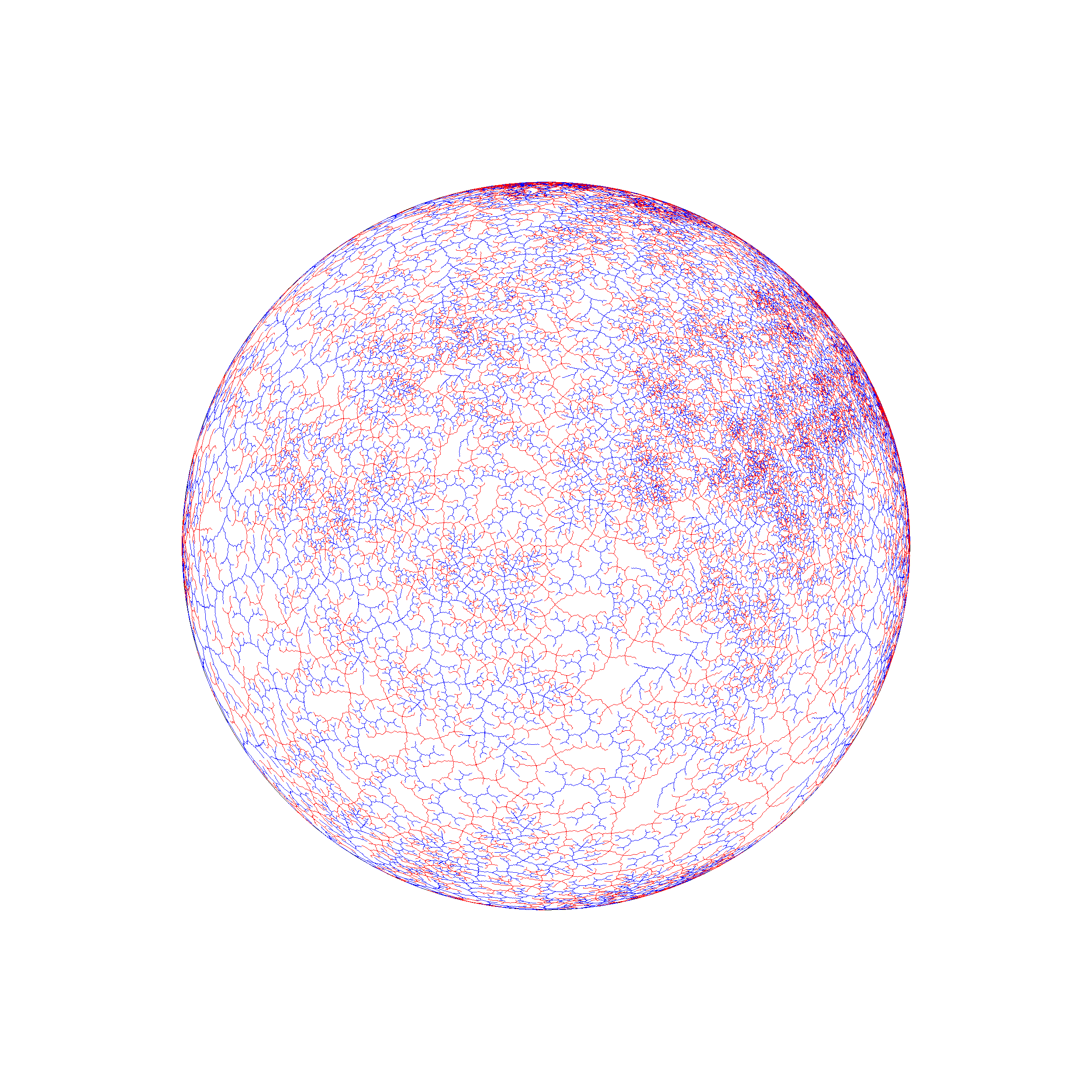}
\caption{A zipper associated to the $(0,2)$ orbifold filling on the Figure 8 knot complement.}
\label{zipper_figure}
\end{figure}

\section{Introduction}

Many well-studied dynamical structures on hyperbolic 3-manifolds (taut foliations, tight essential
laminations, quasigeodesic or pseudo-Anosov flows etc.) are tamed and mastered by a
{\em universal circle} --- a circle with a faithful $\pi_1$-action preserving a pair of
laminations. Under the best circumstances (and conjecturally always) the laminations are
known to be distinct and transverse, and quotienting
the circle by the leaves of the laminations gives rise to a continuous and $\pi_1$-equivariant 
{\em sphere filling Peano curve} from the universal circle to the sphere at infinity 
of hyperbolic 3-space. 

This map from a circle to a sphere is not an embedding. 
But it can be {\em approximated} by an embedding. Rather
fancifully, one might imagine perturbing the image slightly to an embedding, 
and then using this embedded circle as a guide to `unzip' the sphere into two 
very nearly $\pi_1$-invariant path-connected pieces and then somehow taking a limit
under smaller and smaller perturbations.

In fact, this fantasy can be made completely rigorous! If $M$ is a (closed or finite volume) 
hyperbolic 3-manifold, we may simply define a {\em zipper} for $M$ to be a pair $Z^\pm$ of 
disjoint, nonempty, path-connected, $\pi_1(M)$-invariant subsets of $S^2_\infty$ (it should
be emphasized that these subsets are not closed!) The executive summary of our paper has two parts:
\begin{enumerate}
\item{a zipper gives rise to a canonical universal circle; and}
\item{many dynamical structures on hyperbolic manifolds, including some of the
structures mentioned above (certain taut foliations, quasigeodesic flows) as well
as some new ones (uniform quasimorphisms, uniform $\R$-actions) give rise to zippers.}
\end{enumerate}
Even in cases where universal circles were already known to
exist, our construction is new and at least in some cases, arguably simpler and more direct,
and perhaps can play the role of a unifying principle.
Moreover, some of the new cases where our construction works are more directly
tied to certain structures (symplectic representations, left orders) where a direct
connection to foliations, laminations etc. is conjectured and sought (e.g.\/ in the
so-called $L$-space conjecture) but not yet known to exist. Part of our work should 
therefore be interpreted as shedding light on these conjectures.

\subsection{Zippers and Uniforms}

Once one sees the utility of zippers, one is motivated to try to construct them. 
In \cite{Thurston_circles_I} Thurston defined a {\em uniform} foliation to be a
taut foliation of a 3-manifold such that in the universal cover any two leaves are
a uniformly bounded distance apart. This implies, but is slightly stronger than
the conclusion that the foliation is $\R$-covered, which in turn can be expressed
in terms of connectivity properties.

Other structures on a manifold or its fundamental
group give rise to a stratification of the universal cover into (coarse) subsets;
one might say informally that such structures are {\em uniform} if these strata are
coarsely connected. We study two examples in detail --- stratifications coming from
quasimorphisms, and from left orders (or more precisely, orientation-preserving
actions on $\R$) --- and show in either case that a {\em uniform}
structure gives rise to a zipper. 

\subsection{The $L$-space conjecture}

The {\em $L$-space conjecture} says that for an irreducible rational homology 3-sphere $M$ the
following are equivalent:

\begin{enumerate}
\item{the fundamental group $\pi_1(M)$ is left-orderable;}
\item{$M$ is not an $L$-space (i.e.\/ its Heegaard Floer Homology is not minimal); and}
\item{$M$ admits a co-orientable taut foliation.}
\end{enumerate}

Boyer--Gordon--Watson \cite{Boyer_Gordon_Watson} conjectured the equivalence of (1) and (2);
the equivalence of (2) and (3) was posed as a question by Ozsv\'ath and Szab\'o 
\cite{Ozsvath_Szabo}, and conjectured by Juh\'asz \cite{Juhasz}. For an enlightening 
discussion of this conjecture and an enormous amount of experimental evidence bearing on 
it, see Dunfield \cite{Dunfield}.

Whatever the eventual status of the conjecture it has spurred an enormous amount of effort, 
which at the very least has uncovered evidence of a profound and unexpected 
connection between these three kinds of structures on 3-manifolds. 
One of the purposes of this paper is to make a contribution 
to this conjecture by describing a new structure --- a zipper --- 
with manifest connections to all three
legs of the $L$-space conjecture (the connections to taut foliations and left orders should
be obvious; the connections to Heegaard Floer Homology are more indirect, and depend on
the theory of quasimorphisms as a common tool for the investigation of both hyperbolic
and symplectic geometry). 

\subsection{Statement of Results}

For the convenience of the reader we summarize the content of the paper.

In \S~\ref{section:zippers} we introduce zippers for hyperbolic 3-manifolds and
establish their basic properties. The main theorem of this section is
\begin{universal_circle_theorem}
Let $M$ be a hyperbolic 3-manifold and let $Z^\pm$ be a minimal zipper for $M$. 
Then there is a {\em universal circle} $S^1_\un$ associated to $Z^\pm$; in
other words, there is a circle $S^1_\un$ canonically isomorphic to the end circles
of $Z^+$ and $Z^-$, and a faithful action $\pi_1(M) \to \Homeo(S^1_\un)$ 
leaving invariant a pair of laminations $\Lambda^\pm$. 
\end{universal_circle_theorem}

We should stress that we do {\em not} show in this paper that $\Lambda^+$ and $\Lambda^-$
are different. However we expect that this is true, and what we show is that 
each of the laminations $\Lambda^\pm$ separately is related in a precise way to the 
topology of $Z^\pm$ which is too involved to summarize here; 
see \S~\ref{subsection:invariant_laminations} for a more thorough discussion.

In a future paper we plan to prove a
generalization (`Theorem'~\ref{theorem:P_universal_circle}) of this
theorem to a broader class of objects called P-zippers.

In \S~\ref{section:quasigeodesic} we discuss the connection to quasigeodesic flows.
The main theorem of this section is
\begin{quasigeodesic_theorem}
Let $X$ be a quasigeodesic flow on a hyperbolic 3-manifold $M$. The images of the
endpoint map $e^\pm(P^\pm)$ are a P-zipper, and are equal to a zipper 
if and only if $X$ has no perfect fits.
\end{quasigeodesic_theorem}

In \S~\ref{section:quasimorphisms} we define uniform quasimorphisms. A quasimorphism
on a group is uniform if (roughly speaking) the coarse level sets are 
coarsely connected. The main theorem of this section is
\begin{quasimorphism_theorem}
Let $G$ be a hyperbolic group and suppose that $G$ admits a uniform quasimorphism.
Then there is a zipper $Z^\pm \subset \partial_\infty G$ for $G$.
\end{quasimorphism_theorem}

In \S~\ref{section:orders} we define uniform actions. A faithful orientation
preserving action of a group on $\R$ is uniform if
(roughly speaking) the partial order generated by certain {\em up elements} 
(those that move every point in the positive direction) has upper bounds for pairs on a certain scale. 
The main theorem of this section is
\begin{order_theorem}
Let $G$ be a hyperbolic group and suppose $G$ admits a uniform action on $\R$. 
Then there is a zipper $Z^\pm \subset \partial_\infty G$ for $G$.
\end{order_theorem}

\section{Zippers}\label{section:zippers}

In this section we introduce zippers for hyperbolic 3-manifolds, 
and establish their basic properties.
The main theorem of this section (Theorem~\ref{theorem:universal_circle}) 
is that zippers give rise to universal circles.

\subsection{Definition}

Throughout this section we fix a closed oriented 
hyperbolic 3-manifold $M$ and an identification
of the universal cover $\tilde{M}$ with $\H^3$. 
Let $G$ denote the fundamental group of $M$, so that
$G$ acts on $S^2_\infty$ by M\"obius transformations.

\begin{definition}
A {\em zipper} for $M$ is a pair $Z^\pm$ of subsets of $S^2_\infty$ satisfying
\begin{enumerate}
\item{each of $Z^\pm$ is nonempty and path-connected;}
\item{each of $Z^\pm$ is $G$-invariant; and}
\item{$Z^+$ and $Z^-$ are disjoint.}
\end{enumerate}
\end{definition}

\begin{remark}
One may extend this definition to arbitrary hyperbolic groups $G$ 
with Gromov boundary $\partial_\infty G$ (see \cite{Gromov_hyperbolic} for
definitions). A zipper for $G$ is a pair $Z^\pm$ of nonempty path-connected disjoint
$G$-invariant subsets of $\partial_\infty G$. In the sequel we shall see
some constructions of zippers that make sense for arbitrary hyperbolic groups;
however the applications we pursue in this paper 
are only interesting in the case $\partial_\infty G = S^2_\infty$.
\end{remark}

\subsection{Basic properties}\label{subsection:basic_properties}

Fix a zipper $Z^\pm$.

\begin{lemma}[Dense]
Each of $Z^\pm$ is dense in $S^2_\infty$.
\end{lemma}
\begin{proof}
Since $G$ is cocompact, it acts minimally on $S^2_\infty$. Since each of $Z^\pm$ are
$G$-invariant and nonempty, they are dense. 
\end{proof}

\begin{lemma}[Unique path]\label{lemma:unique_path}
Every two points in $Z^+$ are joined by a unique embedded path (and similarly for $Z^-$).
\end{lemma}
\begin{proof}
Suppose not, so that there are distinct points in $Z^+$
that can be joined by two distinct embedded paths. After eliminating backtracks
we can find an embedded circle $\gamma \subset Z^+ \subset S^2_\infty$ which
(by the Jordan curve theorem) separates $S^2_\infty$ into two open disks $D^\pm$. Since
$Z^-$ is nonempty and $G$-invariant it is dense, and therefore it intersects both
$D^\pm$. But $Z^-$ is path-connected, and any path from $D^+$ to $D^-$ must intersect
$\gamma$ and therefore $Z^+$, contrary to the fact that $Z^\pm$ are disjoint. The
proof follows.
\end{proof}

\begin{definition}[Convex hull]
If $P$ is any subset of $Z^+$ (resp $Z^-$), the {\em convex hull} of $P$ is the union of the
embedded paths in $Z^+$ joining pairs of points in $P$.
\end{definition} 
The convex hull of $P$ is the unique minimal
path-connected subset of $Z^+$ containing $P$. If $P$ is finite, its convex hull is
therefore homeomorphic to a finite simplicial tree, as a subspace of $S^2_\infty$. 
Note that every topological embedding of a finite tree in $S^2$ is tame (i.e.\/
homeomorphic to a PL embedding); a standard reference for this is Newman \cite{Newman},
Chapter VI. 

\begin{definition}[Path topology]
The {\em path topology} on $Z^+$ (resp. $Z^-$) is the weak topology generated by the subspace topology
(as a subset of $S^2_\infty$) on all convex hulls of finite subsets. 
\end{definition}

What we call the `weak topology' here is sometimes also called the `coherent topology'; 
i.e.\/ the topology on a union of distinguished subspaces, in which a set is open if 
and only if its intersection with each of the distinguished subspaces is open.

In the path topology,
$Z^+$ (resp. $Z^-$) is a {\em topological $\R$-tree} (also called a {\em dendrite}) 
and the inclusion map $Z^+ \to S^2_\infty$ is continuous but not a homeomorphism 
onto its image.

\begin{definition}[Hard end]
A {\em hard end} in a topological $\R$-tree is a point which is not a cut point. 
\end{definition}
A hard end is not an end in the usual topological sense, but may be related to it as follows.
If $T$ is a topological $\R$-tree, the subspace of $T$ obtained by removing all the hard
ends is itself a topological $\R$-tree $\inte(T)$ called the {\em interior} of $T$, and the
hard ends of $T$ are a subset of the ends (in the usual sense) of $\inte(T)$. One might
express this (by abuse of notation) by referring to cut points as interior points, and then
$\inte(T)$ is just the set of interior points of $T$. Note that
$\inte(\inte(T))=\inte(T)$.

One may obtain the space of ends of $\inte(T)$ directly from $T$ as follows.
An {\em interior subtree} is a finite subtree $K\subset T$ that
contains no hard ends; equivalently $K$ is equal to the convex hull of finitely many
cut points. 

\begin{definition}[End space]
Let $T$ be a topological $\R$-tree. The {\em interior end space} is the inverse limit
$$\EE(T): = \varprojlim \pi_0(T-K)$$
over all interior subtrees $K$ of $T$.
\end{definition}
Note that $\EE(T)=\EE(\inte(T))$, and if $T=\inte(T)$ then $\EE(T)$ 
is equal to the space of ends of $T$ in the usual sense.

\begin{definition}[Minimal zipper]
A zipper $Z^\pm$ is {\em minimal} if $Z^+$ (resp. $Z^-$) is equal to the convex hull of the
$G$-orbit of any point.
\end{definition}
Note that a minimal subzipper necessarily contains no hard ends.

\begin{proposition}[Minimal zipper]\label{proposition:minimal_zippers_exist}
Any zipper contains a minimal subzipper.
\end{proposition}
\begin{proof}
Note first that neither of $Z^\pm$ can contain a global fixed point, since each is
$G$-invariant and $G$ has no global fixed points on $S^2_\infty$.
Furthermore, since $G=\pi_1(M)$ it is finitely generated by $g_1,\cdots,g_n$ (say). Fix
a cut point $p\in Z^+$ and let $X$ be the convex hull of $p \cup \left(\cup_i g_i(p)\right)$.
Then $X$ is a finite tree, and $GX$ (i.e.\/ the union of $G$-translates of $X$) is path
connected and $G$-invariant, and is therefore a $G$-invariant topological $\R$-tree.

For any $q\in X$ let $Y(q)$ denote the convex hull of the $G$-orbit $Gq$ and let
$X(q)$ denote the intersection $Y(q)\cap X$. Then $Y(q)$ is a $G$-invariant topological
$\R$-tree, and $GX(q)=Y(q)$. Since $X$ and $Y(q)$ are both path-connected, and $Z^+$ is
a topological $\R$-tree, it follows that $X(q)$ is path-connected; i.e.\/ it may be
completed to a finite tree by adding finitely many endpoints. 
Since there are no global fixed points, $Y(q)$ has
interior and therefore so does $X(q)$.

Let $\overline{X(q)}$ be the closure of $X(q)$ in $X$ and suppose $q'\in \overline{X(q)}-X(q)$;
i.e.\/ it is one of the endpoints of the finite tree $\overline{X(q)}$.
We claim that $X(q')$ is contained in $\overline{X(q)}$. To show this, it suffices to
show that for any $g\in G$, if $\sigma$ is the embedded path from $q'$ to $gq'$ then the 
interior of $\sigma$ is in $Y(q)$ (since $Y(q')$ is the union of $G$-translates of $\sigma$). 
Since $q'$ is in the closure of $Y(q)$, there are 
$g_1,g_2\in G$ and an embedded path $\tau$ from $g_1q$ to $g_2q$ so that $\tau$ comes
arbitrarily close to $q'$. But then the convex hull of $g_1q, g_2q, gg_1q, gg_2q$ contains
a path whose endpoints come arbitrarily close to $q'$ and $gq'$, and therefore this
convex hull contains `most' of $\sigma$ (and by taking limits, all of the interior of
$\sigma$ is in $Y(q)$). The claim follows.

We may therefore take any maximal nested sequence of subsets $\overline{X(q)}$ and, 
since these are compact nonempty subsets of a compact Hausdorff space, extract a point
$q'$ in the intersection of them all. By the argument above, for any $q''\in Y(q')$
we have $Y(q'')\subset Y(q')$ and by maximality, $Y(q'')$ is dense in $Y(q')$. Therefore
the interior of $Y(q')$ is minimal.
\end{proof}

From now on we will assume that our zippers are minimal. In particular, our zippers
contain no hard ends.

Observe that in the course of the proof of 
Proposition~\ref{proposition:minimal_zippers_exist}
we showed that a minimal zipper is a 
countable union of finite simplicial trees; in particular it
is $\sigma$-compact. We do not use this fact in the sequel.

\subsection{Landing rays}

\begin{definition}[Landing rays]
Let $Z^\pm$ be a zipper. A proper ray $r\subset Z^+$ (resp. $Z^-$) is a {\em landing ray}
if the image in $S^2_\infty$ may be compactified with a single point to a closed
interval. The compactifying point is called the {\em end} of the landing ray, and
we say $r$ {\em lands} at this point (or just that $r$ {\em lands}).
\end{definition}

The terminology `landing rays' comes from the field of holomorphic dynamics.
We use our knowledge of the dynamics of $G$ on $Z^+$ and on $S^2_\infty$ 
to obtain landing rays:

\begin{lemma}[Invariant ray lands]\label{lemma:invariant_ray_lands}
Let $x\in Z^+$ and $g\in G$, and let $\sigma$ be the unique embedded path in $Z^+$
from $x$ to $g(x)$. Then either $\sigma$ contains a fixed point for $g$, or else
$\cup_{n \ge 0} g^n(\sigma)$ contains a ray which lands at the attracting fixed point for $g$
in $S^2_\infty$.
\end{lemma}
\begin{proof}
Suppose $\sigma$ does not contain a fixed point for $g$.
Let $p$ and $q$ be respectively the attracting and repelling fixed point 
of $g$ in $S^2_\infty$. Since $\sigma$ is disjoint from $p \cup q$, the 
forward iterates $g^n(\sigma)$ exit compact subsets of 
$S^2_\infty - (p\cup q)$ and converge only on $p$. This union is path connected, and after
eliminating backtracks, can be seen to contain a landing ray.
\end{proof}

\subsection{Circular orders}

We now come to one of the key properties of zippers: the existence of a canonical 
circular order on each of the end spaces $\EE(Z^\pm)$. 
Informally, these circular orders arises as follows: if our zipper is
minimal, each of $Z^\pm$ is a countable union of finite simplicial trees, each
of which is tamely embedded in $S^2$, and therefore the edges incident to each
vertex inherit a cyclic order from the orientation on $S^2$. These local cyclic orders
are compatible with each other, and determine a cyclic order on ends in an obvious
way. Let's spell this out.

\begin{lemma}[Circular order]
Let $Z^\pm$ be a zipper. The end spaces $\EE(Z^\pm)$ of $Z^\pm$ each admit a 
$G$-invariant circular order.
\end{lemma}
\begin{proof}
Recall that we defined $\EE(Z^+)$ (say) to be the inverse limit of $\pi_0(Z^+-K)$
over finite interior subtrees $K$ of $Z^+$.

Let $e_1,e_2,e_3$ be ends of $Z^+$ associated to
distinct noncompact components of $Z^+-K$ for some specific $K$, and let
$p_1,p_2,p_3$ be points in these components. Let $D$ be a (closed) disk in $S^2_\infty$
containing $K$ in its interior and with $p_1,p_2,p_3$ in its exterior. There are
unique, disjoint paths in $S^2$ from $K$ to the $p_i$, and these intersect $\partial D$
for the first time in three unique points. The orientation on $S^2_\infty$ induces an orientation
on $\partial D$ and therefore a circular order on these three points, and thereby on the
$e_i$. 

\begin{figure}[htpb]
\centering
\includegraphics[scale=0.3]{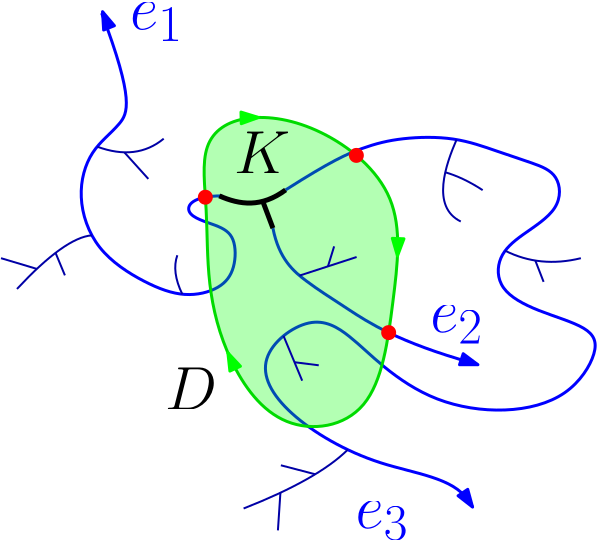}

\includegraphics[scale=0.15]{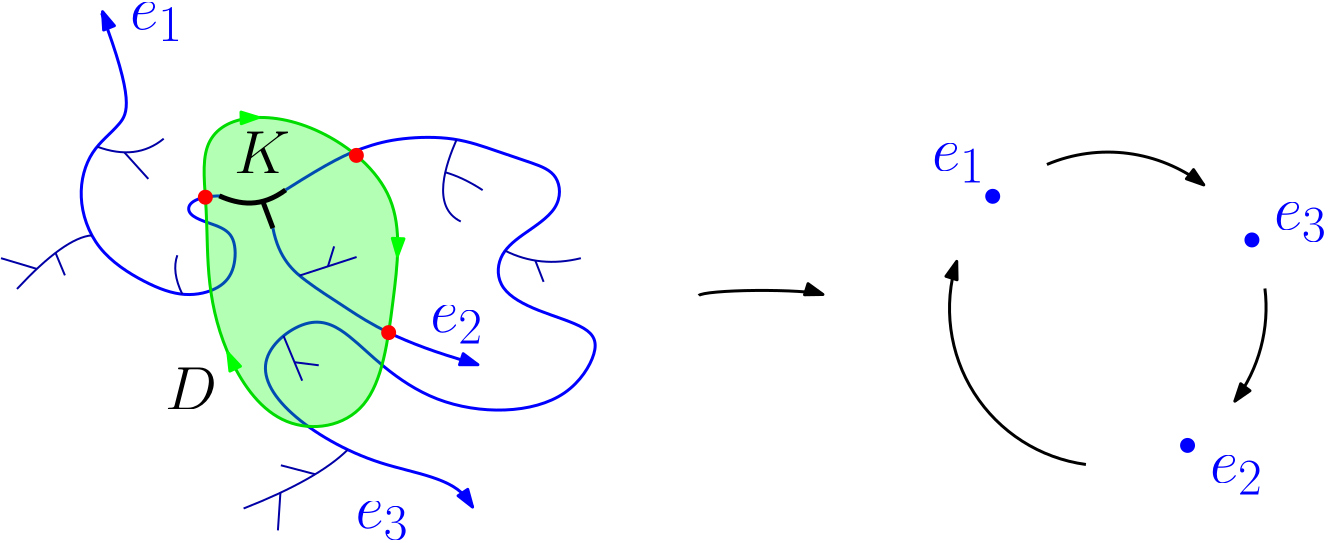}
\caption{There are unique disjoint paths (in blue) from $K$ (in black) to the three ends
$e_1,e_2,e_3$; these intersect $\partial D$ (in green) for the first time in three
distinct points (in red); the cyclic order of these three points in $\partial D$ gives
the circular order on $e_1,e_2,e_3$.}
\label{circular_order_definition}
\end{figure}

This construction is evidently compatible on 4-tuples of ends, and depends only
on the planar structure of the embedding of $Z^+$ in $S^2_\infty$, and is therefore
$G$-invariant.
\end{proof}

Since $Z^+$ is $G$-invariant 
a standard argument implies $\EE(Z^+)$ is either infinite or consists of two points,
and $|\EE(Z^+)|=2$ if and only if $Z^+$ is homeomorphic to an interval.

\begin{lemma}[Infinite ends]\label{lemma:end_stabilizer_cyclic}
The stabilizer of any element of $\EE(Z^\pm)$ is cyclic. Hence in particular
each of $\EE(Z^\pm)$ is infinite.
\end{lemma}
\begin{proof}
Any end with a nontrivial stabilizer lands, by Lemma~\ref{lemma:invariant_ray_lands}.
The stabilizer of this endpoint (and hence of the end) is therefore cyclic, since
$G$ is a Kleinian group without elliptic or parabolic elements. 
\end{proof}

\subsection{Ends and ideal gaps}

The ends of a topological space have a natural inverse limit topology, but it makes
more sense for us to consider a different topology on $\EE(Z^\pm)$, 
namely the (circular) order topology. In this
topology the spaces of ends $\EE(Z^\pm)$ are not necessarily compact, but they 
may be compactified by adding a limit point whenever there is an infinite sequence 
of nested pairs of closed intervals whose intersection in $\EE(Z^\pm)$ is empty.
Such a limit point is represented by geometric data living in $Z^\pm$ that we call an
{\em ideal gap}. Let's explain the construction.

If $e^L,e^R \in \EE(Z^+)$ (say) the closed interval $[e^L,e^R]$ denotes the set of ends
$e\in \EE(Z^+)$ for which $e^L,e,e^R$ are positively cyclically ordered, together
with $e^L,e^R$ themselves. Suppose there is an infinite nested sequence
$$\cdots \subset [e_{i+1}^L,e_{i+1}^R] \subset [e_i^L,e_i^R] \subset \cdots$$
whose intersection is empty. Since the intersection is empty, we may assume without
loss of generality that the $e_i^L$ and $e_i^R$ are all distinct.
The convex hull of each pair $e_i^L,e_i^R$ is a proper
$\R$ in $Z^+$ that we denote $\ell_i$; we orient it from $e_i^L$ to $e_i^R$.

\begin{lemma}[Nonempty intersection]\label{lemma:nested_nonempty_intersection}
Suppose $Z^+$ has no hard ends (for instance, suppose $Z^+$ is minimal) and
suppose further that there is an infinite nested sequence
$$\cdots \subset [e_{i+1}^L,e_{i+1}^R] \subset [e_i^L,e_i^R] \subset \cdots$$
whose intersection is empty. Then with notation as above there is 
an $i$ so that the intersection $X:=\cap_{j\ge i} \ell_j$
is either a point or a compact, nonempty interval.
\end{lemma}
\begin{proof}
Choose a point $p_0 \in \ell_0$ and for each $i$ let $J_i$ be the shortest embedded
interval in $Z^+$ from $p_0$ to some unique nearest $p_i \in \ell_i$ (it is possible
that $p_i = p_0$ and $J_i$ is degenerate). Then because the intervals $[e_i^L,e_i^R]$ are
nested, we have $J_1 \subset J_2 \subset \cdots$ and we may therefore form the union
$J:=\cup J_i$. 

If $J$ were a properly embedded ray, this ray would determine an end which by construction
lies in the intersection of all the $[e_i^L,e_i^R]$, which by hypothesis is empty. Thus
either $J$ is, or may be compactified in $Z^+$ to a compact interval, 
with endpoints $p_0$ and $p_\infty$. But $Z^+$ has no hard ends, so
$p_\infty$ is a cut point, and therefore it must separate some pair of ends. The
only possibility is that it separates $e_i^L$ from $e_i^R$ for all
sufficiently large $i$, which is to say it must lie on all $\ell_i$ for $i$ sufficiently
large. The intersection of any two $\ell_i,\ell_j$ if nonempty is a compact interval, 
and therefore $\cap_{j\ge i} \ell_j$ is compact and connected and, since it contains
$p_\infty$ it is nonempty. This proves the lemma.
\end{proof}

\begin{definition}[End nested; Ideal gap]
A sequence of proper lines $\ell_i$ in $Z^\pm$ which are the convex hulls of the endpoints
of a sequence of nested intervals $[e_i^L,e_i^R]$ in $\EE(Z^\pm)$ is called {\em end nested}.

A point $p\in Z^\pm$ together with the data of an 
infinite sequence of end nested proper lines $\ell_i$ 
with $p \in \cap \ell_i$ 
(as in the statement of Lemma~\ref{lemma:nested_nonempty_intersection}) is called a {\em ideal gap}. 
Two ideal gaps are {\em equivalent} if they determine the same point in the circular
completion of $\EE(Z^\pm)$.
\end{definition}

\begin{figure}[htpb]
\centering
\includegraphics[scale=0.3]{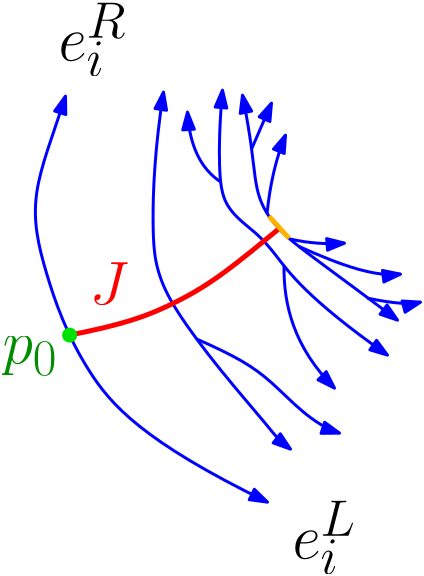}

\includegraphics[scale=0.15]{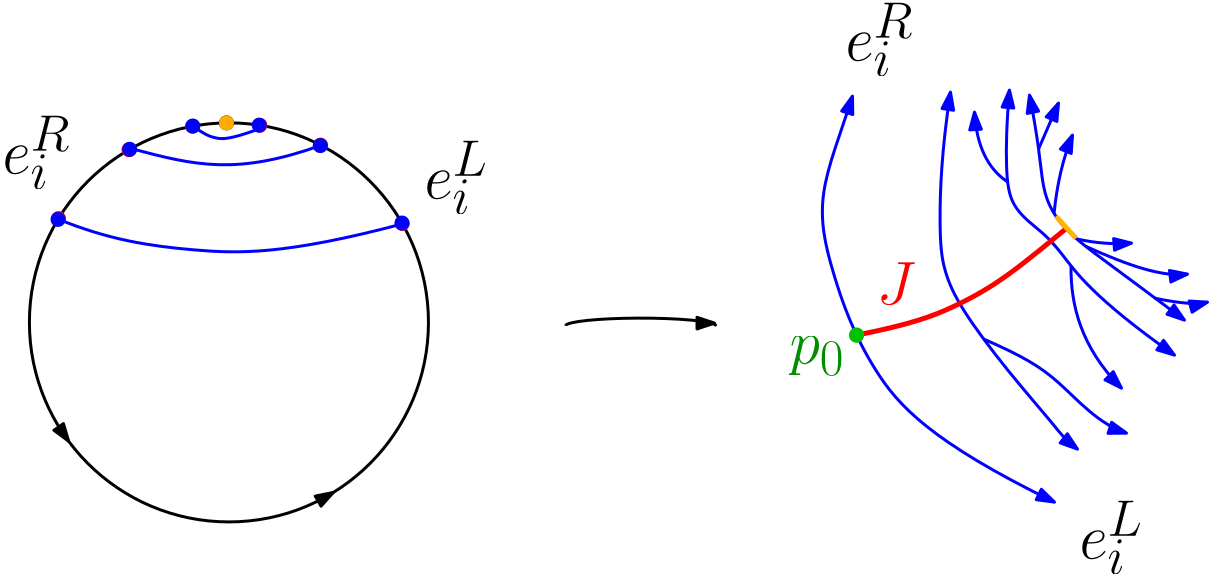}
\caption{$J$ (in red) is a path from $p_0$ (in green) to $X$ (in orange) 
representing the equivalence class of an ideal gap.}
\label{path_to_ideal_gap}
\end{figure}

It is evident from the proof of Lemma~\ref{lemma:nested_nonempty_intersection} that
each equivalence class of ideal gap is of the form $\cup_i \cap_{j\ge i} \ell_i$ for an
infinite end nested sequence of proper lines as above, 
and is either a single point or a connected interval in $Z^\pm$. 

\begin{remark}
It might seem that points representing ideal gaps are rather special in $Z^\pm$, but
this is not so, and in fact it is reasonable to conjecture that {\em every} point in
$Z^\pm$ represents at least two (and possibly infinitely many) ideal gaps. 
Here is the idea. Choose any $p\in Z^+$ minimal (say). Since $p$ is a cut point 
and $Z^+$ is a countable union of segments, there is a
decomposition into countably many components $Z^+-p = \coprod X_i$ each of which
contains an interval of ends $I_i \subset \EE(Z^+)$. Any time there are two consecutive
intervals $I_i, I_j$ so that $I_i$ does not contain a rightmost point and $I_j$ does
not contain a leftmost point, the associated ideal gap is represented by $p$.

We shall return to this construction in \S~\ref{subsection:invariant_laminations}.
\end{remark}

\subsection{End circles}

\begin{proposition}[End circles exist]\label{proposition:end_circles_exist}
Let $Z^\pm$ be a minimal zipper, and let $\overline{\EE}(Z^\pm)$ 
denote the order completions of $\EE(Z^\pm)$. Then
\begin{enumerate}
\item{$\overline{\EE}(Z^\pm)$ are obtained from $\EE(Z^\pm)$ by adding the set of
ideal gaps;}
\item{each space $\overline{\EE}(Z^\pm)$ is separable with countably many (non-ideal) gaps; and}
\item{each space $\overline{\EE}(Z^\pm)$ has no isolated points.}
\end{enumerate}
Consequently by collapsing (non-ideal) gaps we obtain monotone $G$-equivariant surjective maps 
$\overline{\EE}(Z^\pm) \to S^1(Z^\pm)$ whose fibers have cardinality at most $2$, and where the
$S^1(Z^\pm)$ are circles with a natural $G$-action, called the {\em end circles} 
associated to the zipper $Z^\pm$.
\end{proposition}
\begin{proof}
Lemma~\ref{lemma:nested_nonempty_intersection} already implies that the order completion
of $\EE(Z^\pm)$ is obtained by adding ideal gaps. This proves the first bullet.

For a minimal zipper, the action of $G$ on
$\overline{\EE}(Z^\pm)$ is minimal (i.e.\/ every orbit is dense); for otherwise,
the convex hull of the closure of some orbit would be a proper subzipper. It
follows that each of $\overline{\EE}(Z^\pm)$ 
is separable in its order topology. 

A (non-ideal) gap in $\overline{\EE}(Z^\pm)$ is
a pair of distinct points $e_L,e_R$ for which $[e_L,e_R]$ has no interior points; here
each of $e_L,e_R$ is either an end or an ideal gap. The convex hull $\ell$ of $e_L \cup e_R$
in $Z^\pm$ is a nontrivial interval, which might have endpoints or not (it is easy to
see at least one end of $\ell$ is a proper ray, since otherwise $\ell$ would 
represent the equivalence class of a single ideal gap). Furthermore, $\ell$ is
canonically oriented, by the (circular) order. It follows from the definition of a (non-ideal) gap
that if $\ell$ and 
$\ell'$ are intervals associated to distinct gaps, then on any subinterval of the
intersection $\ell \cap \ell'$ the associated orientations {\em disagree}. In particular,
any compact oriented interval $I \subset Z$ embeds in at most one such $\ell$ in an
orientation-preserving way.

Now, any minimal zipper is the $G$ orbit of some compact finite subtree $X$, and we may 
find a countable collection of open oriented subintervals $I_i$ of $X$ so that every embedded 
nontrivial path in $Z$ contains some $G$-translate of some $I_i$. This proves
that the set of (non-ideal) gaps of $\overline{\EE}(Z^\pm)$ is countable.

If $\overline{\EE}(Z^\pm)$ had isolated points, the convex hull of the derived set would
be a proper subzipper (the derived set is nonempty because $\EE(Z^\pm)$ are infinite
by Lemma~\ref{lemma:end_stabilizer_cyclic}). 

Finally, apply \cite{Frankel_circle}, Construction~5.7.
\end{proof}

\begin{remark}
It seems entirely possible that for a minimal zipper, $\overline{\EE}(Z^\pm)$ has no (non-ideal)
gaps at all, so that it is already a circle.
\end{remark}

\begin{proposition}[End circle minimal]\label{proposition:end_circle_minimal}
Let $Z^\pm$ be minimal. Then the action of $G$ on each end circle $S^1(Z^\pm)$ is minimal.
Furthermore point stabilizers are cyclic (and consequently the action is faithful).
\end{proposition}
\begin{proof}
Each point in the end circle $S^1(Z^\pm)$ is either the image of an end or 
an ideal gap in the circular order completion of $\EE(Z^\pm)$. In either case
if the orbit of the point is not dense we can take the convex hull of the associated
end or a point representing the ideal gap and obtain a smaller subzipper, 
contradicting minimality.

Any point in the end circle is the image of a compact subinterval of $\overline{\EE}(Z^\pm)$ 
(note that this compact subinterval consists of either one or
two points) and the stabilizer of the point in the end circle must stabilize the 
endpoints of the corresponding subinterval. These
endpoints are either ends of $Z^\pm$ or ideal gaps represented by points or connected
intervals in $Z^\pm$. The stabilizer of an end is cyclic, by 
Lemma~\ref{lemma:end_stabilizer_cyclic}. The stabilizer of an interval in $Z^\pm$ is
similarly cyclic; for, if this stabilizer is nontrivial, the ends of the interval
must land for a similar reason at points fixed by the entire stabilizer.
\end{proof}

\subsection{Bridges and the universal circle}

The {\em universal circle} $S^1_\un$ is a minimal circle with a $G$ action that collates
$S^1(Z^+)$ and $S^1(Z^-)$ into a single circle. To construct it we must compare the ends
of $Z^+$ and $Z^-$. This is accomplished by means of {\em bridges}.

\begin{definition}[Bridges]
A {\em bridge} is an embedding $\gamma:[0,1] \to S^2_\infty$ such that
\begin{enumerate}
\item{either $\gamma([0,1/2))\subset Z^-$, $\gamma((1/2,1])\subset Z^+$ and
$\gamma(1/2)$ is in $S^2_\infty - (Z^-\cup Z^+)$; or}
\item{$\gamma([0,1))\subset Z^-$ and $\gamma(1)\in Z^+$ or $\gamma([0,1))\subset Z^+$ and
$\gamma(1)\in Z^-$}
\end{enumerate}
where in either case the restrictions to $Z^\pm$ are continuous in the path topology.
The first kind of bridge is said to be of {\em type 1} and the second of {\em type 2}.
\end{definition}
Another way to say this is that a type 1 bridge is a pair of landing rays in $Z^\pm$ that
both land on a common point in $S^2_\infty - (Z^-\cup Z^+)$, and a type 2 bridge is a
landing ray in one of $Z^\pm$ that lands on a point in the other one.

\begin{proposition}[Bridges exist]\label{proposition:bridges}
A bridge exists. 
\end{proposition}

This is proved by an analysis of cases. Every element $g\in G$ is a M\"obius transformation,
and fixes exactly two distinct points in $S^2_\infty$. It follows that each
$g\in G$ fixes at most two points in each of $Z^\pm$ and therefore either some
$g$ fixes zero points in one of $Z^\pm$ or every nontrivial $g$ fixes exactly one point
in each of $Z^\pm$.

\begin{lemma}[fixed point free has tame axis]\label{lemma:axis}
Let $g$ act on $Z^+$ without any fixed points. Then $g$ has a proper axis $\ell \subset Z^+$
which is homeomorphic to $\R$, and whose ends are both landing rays that land at the
fixed points of $g$ in $S^2_\infty$.
\end{lemma}
\begin{proof}
Pick any $p\in Z^+$ and let $\sigma$ be the unique embedded path from $p$ to $g(p)$.
Since $\sigma$ is compact and disjoint from the fixed points of $g$ in $S^2_\infty$,
some finite power of $g$ takes $\sigma$ off itself. It follows that the union
$\cup_{n \in \Z} g^n(\sigma)$ contains a unique embedded minimal $g$-invariant axis 
which is $\ell$. The ends of $\ell$ are $g$-invariant, and therefore land by 
Lemma~\ref{lemma:invariant_ray_lands}.
\end{proof}

\begin{lemma}\label{lemma:ray}
Let $G$ act on $Z^+$ in such a way that every nontrivial element has exactly one fixed point.
Then there is some nontrivial $g\in G$ fixing $p^+\in Z^+$ and leaving invariant some 
proper ray $r$ in $Z^+$ starting at $p^+$ that lands on the other fixed point of $g$ in
$S^2_\infty$.
\end{lemma}
\begin{proof}
Suppose the lemma is false. We claim in this case that every nontrivial element $g \in G$ 
fixes exactly one point $p(g)\in Z^+$ and freely permutes the components of $Z^+-p(g)$.
For otherwise there is some component $X \subset Z^+-p(g)$ which is fixed by $g$. In this
case, pick $q\in X$ and let $\sigma$ be an embedded path from $p(g)$ to $q$. Then
$\sigma \cap g(\sigma)$ is nonempty (or else $q$ and $g(q)$ would be in different components)
and therefore contains some initial common interval $\tau$ so that after possibly replacing
$g$ by $g^{-1}$ we may assume $g(\tau) \subset \tau$. But then the union $r:=\cup_{n\le 0} g^n(\tau)$
is a $g$-invariant ray, which consequently lands on a fixed point by 
Lemma~\ref{lemma:invariant_ray_lands}, and the lemma would be true.

\begin{figure}[htpb]
\labellist
\small\hair 2pt
\pinlabel $p$ at 150 35
\pinlabel $q$ at 250 35
\pinlabel $\sigma$ at 200 36
\pinlabel $g(q)$ at 50 28
\pinlabel $h(p)$ at 350 28
\pinlabel $g(\sigma)$ at 100 33
\pinlabel $h(\sigma)$ at 300 33
\endlabellist
\centering
\includegraphics[scale=1]{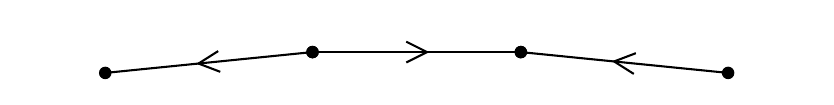}
\caption{$gh^{-1}$ has an axis.}
\label{axis_exists}
\end{figure}

So let's suppose every nontrivial element of $G$ fixes exactly one point in $Z^+$ and
freely permutes the components of the complement of this point.
But now let $g,h$ be elements with distinct fixed points $p,q$ that cobound an embedded 
interval $\sigma$. Evidently $g(\sigma)$ and $h(\sigma)$ are intervals disjoint from each
other and sharing endpoints only with $\sigma$. Then $h(\sigma)$ and $g(\sigma)$ are disjoint
and both contained in an embedded oriented interval in an orientation compatible way
(see Figure~\ref{axis_exists}) and therefore $gh^{-1}$ has an axis composed of 
translates of this interval (in particular it is fixed-point free on $Z^+$, 
contrary to hypothesis). This contradicts the claim, and therefore proves the lemma.
\end{proof}

We now give the proof of Proposition~\ref{proposition:bridges}.
\begin{proof}
If some $g$ acts on $Z^+$ (say) without any fixed points, we may construct a bridge
as follows. If $g$ fixes some point $p\in Z^-$, by Lemma~\ref{lemma:axis}
we may build a type 2 bridge from
part of an axis for $g$ in $Z^+$ together with $p\in Z^-$.
Otherwise $g$ has axes in both $Z^+$ and $Z^-$, and we may construct a type 1 bridge from
part of each of these together with a fixed point for $g$ in $S^2_\infty$.

Otherwise we may assume every nontrivial element of $G$ fixes exactly one point in
each of $Z^\pm$. In this case by Lemma~\ref{lemma:ray} we may build a type 2 bridge from a $g$-invariant
ray in $Z^+$ together with the other fixed point of $g$ in $Z^-$ (or vice versa).
\end{proof}

A bridge determines a correspondence between certain points in $S^1(Z^\pm)$. 
A type 1 bridge consists of two proper rays in $Z^\pm$; these rays determine elements of
$\EE(Z^\pm)$ and therefore points in $S^1(Z^\pm)$.
A type 2 bridge consists of a ray $r$ in $Z^+$ (say) and a point $p$ in $Z^-$. The ray in $Z^+$
determines an element of $\EE(Z^+)$ and therefore a point in $S^1(Z^+)$. This same ray
lands at a point $p \in Z^-$ but is disjoint from the components of $Z^--p$; if we attach
$r$ to $Z^-$ the result $Z^-\cup r$ is a topological $\R$-tree with a unique hard end
at the initial point of $r$, and this end has a well-defined place in the circularly
ordered set $\EE(Z^-\cup r)$ which refines the circular order on $\EE(Z^-)$. In particular,
it determines a point in the order completion of $\EE(Z^-)$
and therefore a point in $S^1(Z^-)$.
We refer to this pair of points in $S^1(Z^\pm)$ in either case as the {\em impression} of the bridge.

The collection of all bridges determines in this way a $G$-equivariant correspondence between
nonempty subsets of $S^1(Z^\pm)$ that we call the {\em bridge correspondence}.

\begin{lemma}[Bridge correspondence reverses]\label{lemma:order_reverses}
The bridge correspondence reverses the circular orders on $S^1(Z^\pm)$.
\end{lemma}
\begin{proof}
Let $p^\pm, q^\pm, r^\pm \in S^1(Z^\pm)$ be paired by the bridge correspondence. Each of these points
either lies in $Z^\pm$ or is an end of $Z^\pm$; thus each triple spans a convex hull in 
$Z^\pm$ which is a (possibly degenerate) tripod that may be compactified (if necessary) by adding
endpoints in $S^2_\infty$. Let $X^\pm$ be these tripods. The union of these two tripods is a
theta graph in $S^2_\infty$, and the circular orders on $p^+,q^+,r^+$ resp. $p^-,q^-,r^-$
agree with the cyclic orderings of the edges of the theta graph at one resp. the other vertex; 
in particular they disagree with each other.
\end{proof}

\begin{proposition}[Bridge homeomorphism]
Let $Z^\pm$ be a minimal zipper. The bridge correspondence determines a 
$G$-equivariant orientation-reversing homeomorphism between $S^1(Z^\pm)$.
\end{proposition}
\begin{proof}
Since $Z^\pm$ are minimal, the $G$ action on each circle $S^1(Z^\pm)$ is minimal
by Proposition~\ref{proposition:end_circle_minimal}. We may therefore choose a
bridge together with its $G$-orbit and observe that the set of impressions of this family
of bridges is dense in either circle. Furthermore the stabilizer of the impression
in each circle is, by construction, the same cyclic subgroup of $G$.
Since by Lemma~\ref{lemma:order_reverses} the bridge correspondence reverses
circular orders, and it is a bijection between dense subsets, 
it extends uniquely to an orientation-reversing homeomorphism.
\end{proof}

We may therefore identify both $S^1(Z^\pm)$ with a single $G$-circle, called the
{\em universal circle} $S^1_\un$. Note that although there is no canonical orientation on
$S^1_\un$, it is nevertheless true that the action of $G$ is orientation-preserving.

\begin{remark}\label{remark:Kim_1}
It should be evident that fixed points of group elements play a key role in 
understanding the geometry and topology of zippers. It would be worthwhile to
develop a clearer understanding of this relationship. In a previous version of this
paper we explicitly asked the following questions, implicitly about a minimal zipper $Z^\pm$:
\begin{enumerate}
\item{is it possible that some (nontrivial) element $g\in G$ has both fixed points in $Z^+$?} 
\item{can there be elements $g\in G$ with one fixed point in $Z^+$ and the other fixed point in neither of $Z^\pm$?}
\end{enumerate}
The first question has been answered by KyeongRo Kim with a remarkably beautiful argument
\cite{Kim}, who shows that this possibility
does {\em not} occur: for a minimal zipper, no nontrivial element $g\in G$ can have two fixed points
in $Z^+$ (or $Z^-$). In the same preprint he puts strong constraints on the possibilities for
the second question, making a negative answer seem very likely there too. For example,
Kim is able to answer a weak version of the second question by showing that there is
{\em some} nontrivial element $g\in G$ which has one fixed point in $Z^+$ and the other
in $Z^-$.
\end{remark}

\subsection{Invariant laminations}\label{subsection:invariant_laminations}

We now explain how to obtain a pair of invariant laminations of $S^1_\un$. For
definitions and an introduction to the theory of laminations of the circle,
see e.g.\/ \cite{Calegari_foliations_book}, \S~2.1.

For each $p\in Z^+$ we get a collection of components $X$ of $Z^+-p$, and each such component
determines a nonempty subset $\EE(X)$ of $\EE(Z^+)$. The closure of each such subset in 
$S^1(Z^+)$ is a nontrivial interval since no end is isolated 
(by Proposition~\ref{proposition:end_circles_exist}, bullet (3)), 
and the endpoints of each such interval are an unordered pair 
$\ell$ of distinct points in $S^1(Z^+)$. By construction, no two pairs $\ell$, $\ell'$ associated to any
two $p,q\in Z^+$ may link in $S^1(Z^+)$ (though they may share an endpoint in common). The
closure of the union of such pairs over all $p\in Z^+$ and all components of $Z^+-p$ is therefore
a $G$-invariant lamination. Using the canonical identification between $S^1(Z^+)$ and $S^1_\un$ we
obtain in this way a $G$-invariant lamination of $S^1_\un$, and denote it $\Lambda^+$. In a similar
way we obtain $\Lambda^-$ from $Z^-$.

\begin{example}[Non-branch points and simple leaves]\label{example:2_valent_point}
A minimal zipper is a countable increasing union of finite trees; it follows that there are
only countably many branch points in each of $Z^\pm$ (i.e.\/ points whose complements have
more than 2 components), and non-branch points are dense. If $p \in Z^+$ is a non-branch
point then $Z^+ - p$ has exactly two components (by definition) we get a unique leaf
$\ell_p \in \Lambda^+$ associated to such a point $p$, and similarly for a non-branch point
in $Z^-$. Call such leaves {\em simple leaves}.
\end{example}

\begin{remark}\label{remark:Kim_2}
If $Z^+$ contains an interval $I$ of non-branch points, then $\ell_p = \ell_q$ 
for any two points $p,q\in I$. We suspect in fact that this can never happen, and this
should follow from the techniques of Kim \cite{Kim}. However, since we have no 
applications in mind, we do not pursue this here.
\end{remark}

\begin{proposition}[Simple leaves]\label{proposition:simple_leaves}
Simple leaves are dense in $\Lambda^\pm$. Furthermore, a leaf of $\Lambda^\pm$ is simple
if and only if it is a limit of leaves from either side. 
\end{proposition}
\begin{proof}
If $p \in Z^+$ is arbitrary, and $X$ is a component of $Z^+-p$, then
for any sequence of non-branch points $p_i \in X$ converging to $p$ the associated simple leaves
$\ell_{p_i}$ converge to the leaf of $\Lambda^+$ associated to the component $X$. This proves
that simple leaves are dense.

If $\ell_p \in \Lambda^+$ is a simple leaf associated to a non-branch point $p\in Z^+$,
the argument above (applied to each of the two components of $Z^+-p$) implies $\ell_p$ is
a limit of leaves from either side.

Conversely, suppose a leaf $\ell\in \Lambda^+$ is a limit from either side. Then it is a limit of simple
leaves from either side, say $\ell_{p_i}$ from one side and $\ell_{q_i}$ from the other.
Evidently the compact intervals $I_i$ in $Z^+$ from $p_i$ to $q_i$ are nested, and 
any point $p \in I:=\cap_i I_i$ is a non-branch point with $\ell_p = \ell$.
\end{proof}

To say anything more precise about the laminations $\Lambda^\pm$ one would need to know
more about the dynamics of $G$ on $Z^\pm$ and on $S^1_\un$. In particular, we do not know the answers to
the following questions:
\begin{enumerate}
\item{Let $Z^\pm$ be a minimal zipper. Are $\Lambda^\pm$ minimal? i.e.\/ is the $G$-orbit of
every leaf dense (in either lamination)?}
\item{Are $\Lambda^\pm$ transverse? Is it possible that they are equal?}
\item{Do $\Lambda^\pm$ have perfect fits? (see \S~\ref{subsection:perfect_fits}).}
\end{enumerate}

\begin{remark}
For the interested reader we briefly give an informal sketch of an argument addressing
the second question, and depending on the results and methods of Kim \cite{Kim} alluded
to earlier.

Kim shows that there is always some nontrivial element $g\in G$ with exactly one
fixed point $p^\pm$ in each of $Z^\pm$. Suppose it is further true that each of $p^\pm$
is the endpoint of a type 2 bridge which we may assume has the other fixed point as
its initial point. That is, there should be proper rays $r^\pm$ in $Z^\pm$ with
initial point $p^\pm$ which limit to $p^\mp$ respectively. Assuming a positive
answer to the question raised in Remark~\ref{remark:Kim_2} every proper ray $r$ in $Z^\pm$
gives rise to a {\em rainbow} --- i.e.\/ a sequence of nested leaves in $\Lambda^\pm$
with distinct endpoints converging to the image of $r$ in $S^1_\un$. In this case,
the limit points of the rainbows associated to $r^\pm$ are the impressions of the
associated bridges; call them $x^\pm \in S^1_\un$. 
One may observe that the limit of the images of leaves in each rainbow
under suitable positive or negative powers of $g$ are nontrivial leaves of $\Lambda^\pm$
whose endpoints in $S^1_\un$ are invariant under $g$ and distinct from $x^\pm$ respectively; 
in particular, $g$ must have at least 3 fixed points in $S^1_\un$.
On the other hand, $x^+$ is not a limit point of a negative rainbow, and similarly
for $x^-$. To see this, suppose not and choose any sequence of leaves of $\Lambda^-$ that 
nest down to $x^+$ and are associated to components of $Z^- - p_i$ for suitable points 
$p_i$. It is easy to see such components must simultaneously approach $x^+$ from `either 
side' of $r^+$. But once a component intersects the interior of the disk with boundary
$r^+ \cup r^-$ it must either intersect $r^-$ (in which case the associated leaf has an
endpoint bounded away from $x^+$) or it cannot intersect the exterior, thus giving a 
contradiction in either case. It follows that $\Lambda^+$ and $\Lambda^-$ are not equal.
\end{remark}

There is an extensive literature on the theory of so-called laminar groups (i.e.\/ groups
that act on a circle in such a way as to preserve a lamination) and there is a range of
subtle phenomena that can occur. For a recent survey article on this subject we suggest
\cite{Baik_Kim}.

\subsection{Main Theorem}

Putting together the results from the previous subsections gives the main theorem of this
section:

\begin{theorem}[Universal circle]\label{theorem:universal_circle}
Let $M$ be a hyperbolic 3-manifold and let $Z^\pm$ be a minimal zipper for $M$. 
Then there is a {\em universal circle} $S^1_\un$ associated to $Z^\pm$; in
other words, there is a circle $S^1_\un$ canonically isomorphic to the end circles
of $Z^+$ and $Z^-$, and a faithful action $\pi_1(M) \to \Homeo(S^1_\un)$ 
leaving invariant a pair of laminations $\Lambda^\pm$. 
\end{theorem}

\begin{remark}
If $Z^\pm$ is not minimal then a priori different choices of minimal subzippers might
give rise to non-isomorphic universal circles. We do not know if this possibility 
actually occurs.
\end{remark}

\subsection{P-zippers}

Certain natural constructions give rise, not to zippers, but to a 
generalization called P-zippers:

\begin{definition}[P-zippers]\label{definition:P_zipper}
A {\em P-zipper} is a pair $Z^\pm$ of subsets of $S^2_\infty$ satisfying
\begin{enumerate}
\item{each of $Z^\pm$ is nonempty and path-connected;}
\item{each of $Z^\pm$ is $G$-invariant;}
\item{there are path connected spaces $Y^\pm$ with $G$ actions and
$G$-equivariant continuous surjective maps $e^\pm:Y^\pm \to Z^\pm$; and}
\item{for every pair of continuous paths $a\subset Y^+$ and $b\subset Y^-$ the
images $\alpha:=e^+(a)$ and $\beta:=e^-(b)$ do not cross in $S^2_\infty$.}
\end{enumerate}
\end{definition}
The last condition means that there is no closed disk $D \subset S^2_\infty$ and 
paths $\alpha$, $\beta$ as above properly contained in 
$D$ so that $\alpha \cap \partial D$ and $\beta \cap \partial D$ are 
disjoint linked $S^0$s in $S^1$.
The `P' in P-zippers may be read as standing equally for `Peano', `parametric', 
or `perfect fits'. Since each of $Z^\pm$ is $G$-invariant, they are dense.
Let's refer to a path $I \to Z^\pm$ that factors through $I \to Y^\pm \to Z^\pm$ as
a {\em liftable path}. Thus the last condition just says that liftable paths cannot
cross.

\begin{example}
Every zipper is a P-zipper.
\end{example}

\begin{example}
P-Zippers do not have to be disjoint, or even distinct! Suppose
$Y^+=Y^-=S^1$ and $e^+=e^-$ is a Cannon--Thurston map associated to a fibering of
a 3-manifold over the circle. In this case the maps $e^\pm$ may be perturbed to be
embeddings with disjoint image; in particular, no liftable paths can cross, since
crossing (which forces the images to intersect) is stable under perturbation.
\end{example}

Thus the idea of a P-zipper is a structure that interpolates between 
(ordinary) zippers and Peano circles, and generalizes both.

In a future paper (work in progress) we expect to prove the following 
parallel to Theorem~\ref{theorem:universal_circle}:

\begin{theoremm}[Universal circle]\label{theorem:P_universal_circle}
Let $Z^\pm$ be a P-zipper for $M$. Then there is a {\em universal circle} $S^1_\un$,
and a faithful action $\pi_1(M) \to \Homeo(S^1_\un)$ leaving invariant a 
pair of laminations $\Lambda^\pm$. 
\end{theoremm}

Roughly speaking the idea of the proof is as follows. First of all, without loss of
generality we may assume that each of $Y^\pm$ is the Cayley graph of $\pi_1(M)$ 
with respect to any finite generating set. Because liftable paths cannot cross,
the restriction of $e^\pm$ to any embedded loop in $Y^\pm$ must factor through
a quotient topological $\R$-tree, and by approximating the maps on these trees by embeddings
one may obtain a local circular ordering of the complementary components at every local cut point.
From this one can construct end circles and laminations. Comparing end circles can
be done by a generalization of the method of bridges.

\subsection{Peano curves}

A universal circle $S^1_\un$ embeds as the equator of $S^2$, and the invariant
laminations $\Lambda^\pm$ embed as geodesic laminations on either side. One may obtain
a quotient space by collapsing leaves of $\Lambda^\pm$ and their complementary gaps to
points. Under many circumstances the resulting quotient space is homeomorphic to $S^2$
by a theorem of R. L. Moore \cite{Moore}.

Since the construction involves no choices, the fundamental group $\pi_1(M)$ acts on this
quotient sphere. If one can show that this is a convergence action, it follows by
a theorem of Bowditch \cite{Bowditch} that the action is conjugate to the action of
$\pi_1(M)$ on its Gromov boundary $S^2_\infty$. In this case, one obtains a continuous
equivariant surjection $S^1_\un \to S^2_\infty$. This construction holds for the
universal circles associated to surface bundles by Cannon--Thurston \cite{Cannon_Thurston},
for universal circles associated to certain pseudo-Anosov flows by Fenley \cite{Fenley_flow_circle},
and for universal circles associated to quasigeodesic flows by Frankel 
\cite{Frankel_Cannon_Thurston}.

It seems very plausible to us that the construction should go through for an arbitrary
zipper; thus we make the

\begin{conjecture}
Let $Z^\pm$ be a (P-)zipper for $M$ with universal circle $S^1_\un$.
There is a continuous $\pi_1$-equivariant surjection $S^1_\un \to S^2_\infty$.
\end{conjecture}

\begin{remark}
The theory of CaTherine wheels (\cite{Calegari_Loukidou}, in preparation) aims to 
shed light on this conjecture from the other direction. It starts with a surjective Peano curve
$f:S^1 \to S^2$ satisfying a small list of axioms, and from this one reconstructs
functorially a zipper $Z^\pm$ and a pair of invariant laminations $\Lambda^\pm$.
\end{remark}

\section{Quasigeodesic flows}\label{section:quasigeodesic}

The main theorem of this section (Theorem~\ref{theorem:zippers_from_flows}) 
is that quasigeodesic flows in general give rise to P-zippers, and those
that satisfy a further technical condition --- those without {\em perfect fits} ---
give rise to (honest) zippers.
This theorem depends on papers of Fenley, Frankel 
and the first author (as we explain) and we refer the reader to these papers for
background and many technical details of the theory of quasigeodesic flows.

\subsection{Quasigeodesic flows}

\begin{definition}[Quasigeodesic flow]
Let $M$ be a closed or finite volume complete hyperbolic 3-manifold. 
A {\em quasigeodesic flow} is a nonsingular
vector field $X$ on $M$ so that if $\tilde{X}$ denotes the lift of $X$ to the universal cover
$\tilde{M}=\H^3$, the flowlines of $\tilde{X}$ are (uniformly) quasigeodesic in $\H^3$.
\end{definition}

The flowlines of $\tilde{X}$ are the leaves of a 1-dimensional foliation of $\H^3$. Denote
the leaf space of this foliation by $P$. Each leaf $\ell$ of this foliation is a quasigeodesic,
oriented by $\tilde{X}$, and therefore limits in forward and backward time 
to unique points $e^\pm(\ell)\in S^2_\infty$.
In \cite{Calegari_quasigeodesic} it is shown that $P$ is homeomorphic to $\R^2$, and that
the maps $e^\pm:P \to S^2_\infty$ are continuous.

Each of the maps $e^\pm$ determines a decomposition of $P$ into the connected components
of point preimages. These components are noncompact, and their ends are circularly ordered;
thus one may obtain a universal circle $S^1_\un$ collating them (see \cite{Calegari_quasigeodesic},
Theorem~B).

Frankel showed (\cite{Frankel_Cannon_Thurston}, Continuous Extension Theorem) 
that one may $G$-equivariantly compactify two copies $P^\pm$ of $P$
by a common universal circle $S^1_\un$ to build a {\em universal sphere} $S^2_\un$ in such a way
that the maps $e^\pm:P^\pm \to S^2_\infty$ extend continuously to a common map 
$e:S^2_\un \to S^2_\infty$. Note that the restriction $e:S^1_\un \to S^2_\infty$ is a
$G$-equivariant surjection; thus Frankel's Theorem generalizes the well-known Cannon--Thurston
theorem for surface bundles.

\subsection{Perfect fits and P-zippers}\label{subsection:perfect_fits}

The point preimages of $e^\pm$ determine a relation on $S^1_\un$
that may be encoded as a pair of $G$-invariant laminations $\Lambda^\pm$
which extends the decomposition associated to $e^\pm$ on each of $P^\pm$ (again, see
\cite{Calegari_quasigeodesic}, Theorem~B).

Evidently the images $e^+(P^+)$ and $e^-(P^-)$ are path-connected and $G$-invariant, and
therefore determine a zipper if and only if they are disjoint. This is captured precisely
by a condition known as perfect fits:

\begin{definition}[Perfect fits]
A quasigeodesic flow has {\em perfect fits} if there are leaves of $\Lambda^\pm$ with
the same endpoint in $S^1_\un$.
\end{definition}

The terminology comes from the theory of pseudo-Anosov flows where it was introduced by
Fenley (see e.g.\/ \cite{Fenley_good_geometry}, Definition~4.2) who showed that perfect fits are the
key to understanding the large scale geometric properties of pseudo-Anosov flows,
quasigeodesic or otherwise. Actually, the connection between the two uses of perfect
fits (in the theories of quasigeodesic and pseudo-Anosov flows respectively) 
is even closer than an analogy, since Frankel and Landry have shown (\cite{Frankel_Landry}) that
any quasigeodesic flow may be deformed (keeping the universal circle constant) to a
quasigeodesic pseudo-Anosov flow; thus a quasigeodesic flow has perfect fits if and only
if the associated pseudo-Anosov flow does.

\begin{theorem}[(P)-Zippers from flows]\label{theorem:zippers_from_flows}
Let $X$ be a quasigeodesic flow on a hyperbolic 3-manifold $M$. The images of the
endpoint map $e^\pm(P^\pm)$ are a P-zipper, and are equal to a zipper 
if and only if $X$ has no perfect fits.
\end{theorem}
\begin{proof}
If a quasigeodesic flow $X$ has no perfect fits, the compactification of any two decomposition 
elements of $e^\pm$ are disjoint in $S^1_\un$, and therefore their images in $S^2_\infty$
are disjoint. In other words, $e^+(P^+)$ and $e^-(P^-)$ are disjoint, and are consequently
zippers. Conversely, if $X$ has perfect fits, then directly from the definition and 
Frankel's Theorem the images $e^+(P^+)$ and $e^-(P^-)$ are not disjoint.

Whether $X$ has perfect fits or not, we claim $e^+(P^+)$ and $e^-(P^-)$ are a
P-zipper, where we take $Y^\pm = \H^3$ and by abuse of notation denote the composition
of the projections $\H^3 \to P^\pm$ with $e^\pm$ by $e^\pm:\H^3 \to S^2_\infty$. 
Let $a$ and $b$ be arcs in $\tilde{M}=\H^3$ with
images $\alpha=e^+(a)$ and $\beta=e^-(b)$. We claim that $\alpha$ and $\beta$ cannot 
cross in the sense of Definition~\ref{definition:P_zipper}. Let us suppose not, and
obtain a contradiction.

Let $A:I \times [0,1) \to \H^3$ be a proper embedding so that $A(I\times 0)$ parameterizes the
arc $a$, and each $A(p\times [0,1))$ is the positive flowline of $\tilde{X}$ emanating
from $A(p,0)\in a$. The map $A$ may be compactified to a map 
$\bar{A}:I\times I \to \H^3 \cup S^2_\infty$
for which $A(I\times 1)$ parameterizes $\alpha$.

Likewise, let $B:I\times (0,1] \to \H^3$ be a proper embedding so that $B(I\times 1)$ parameterizes
$b$, each $B(p\times (0,1])$ is the negative flowline of $\tilde{X}$ emanating
from $B(p,1)\in b$, and $B$ may be compactified to $\bar{B}:I\times I \to \H^3 \cup S^2_\infty$
by $B(I\times 0) \to \beta$.

Now let's think about the intersection of (the images of) $A$ and $B$. A positively
oriented ray and a negatively oriented ray contained in flowlines of $\tilde{X}$ may
intersect in (at most) a compact subinterval. Since the image of $A$ resp. $B$ are 
compact families of positive resp. negative rays, it follows that the intersection of
(the images of) $A$ and $B$ is a {\em compact subset} of $\H^3$. Therefore
there are $t_0 \in [0,1)$ and $t_1 = (0,1]$ so that $A:I\times [t_0,1) \to \H^3$ 
and $B:I\times (0,t_1] \to \H^3$ are {\em disjoint}.

But $\bar{A}(I\times 1)=\alpha$ and $\bar{B}(I\times 0)=\beta$ cross essentially in
$S^2_\infty$ (in the sense of Definition~\ref{definition:P_zipper}). 
This means that if we perturb $\alpha$ rel. endpoints to some
arc $\alpha' \subset \bar{A}(I\times [t_0,1]) \subset \H^3 \cup S^2_\infty$,
properly embedded in $\H^3$ except for its endpoints, the resulting arc will nontrivially
link the circle $\bar{B}(\partial (I\times [0,t_1]))$. In particular, $\alpha'$ has nontrivial
intersection number with the interior of the square $\bar{B}(I\times [0,t_1])$, so
that $\alpha' \cap B((0,1)\times (0,t_1))$ is nonempty. But this means the images of
$A$ and $B$ are not disjoint after all; this contradiction shows that $\alpha$ and
$\beta$ do not cross, as claimed.
\end{proof}

\begin{example}[Surface bundles]\label{example:surface_bundles}
Let $M$ be a closed hyperbolic 3-manifold that fibers over the circle with fiber $S$ 
and monodromy $\phi$. Any flow $X$ transverse
to a foliation by fibers is quasigeodesic; in fact, if $\alpha$ is a non-singular closed
1-form on $M$ with $\ker(\alpha)$ tangent to the surface fibers, the integral of $\alpha$
along any flowline of $X$ increases linearly with the parameterization, certifying that
the lifts of the flowlines to $\H^3$ make uniformly linear progress.

Any two (uniform) quasigeodesic rays in $\H^3$ with the same endpoint are a bounded
distance apart. However, if $r$ and $s$ are proper rays contained in flowlines of
$\tilde{X}$, and if $r$ is positively oriented whereas $s$ is negatively oriented, the
integral of $\alpha$ along $r$ diverges to $+\infty$ while the integral of $\alpha$ along
$s$ diverges to $-\infty$, so that these rays must themselves diverge in $\H^3$ and
can't limit to the same point.
It follows that $X$ as above never has perfect fits, and any fibration of $M$ over $S^1$ 
gives rise to a zipper. 

Figure~\ref{zipper_figure} gives an illustration of the zipper
associated to the fibration over the circle with fiber a torus with one orbifold point of
order 2, and monodromy $\left( \begin{smallmatrix} 2 & 1 \\ 1 & 1 \end{smallmatrix} \right)$.
\end{example}

\begin{example}[Noncompact surface bundles]
Let $M$ fiber over the circle with fiber $S$ and monodromy $\phi$ where $S$ is
{\em noncompact} of finite type; i.e.\/ $S$ is homeomorphic to a closed surface with
finitely many points removed. Then Hoffoss \cite{Hoffoss} showed the suspension flow
of the pseudo-Anosov map is uniformly quasigeodesic in the hyperbolic metric. Thus the
flow $X$ gives rise to a P-zipper in this case too.
\end{example}

\begin{example}[Higher dimensions]\label{example:higher_dimensions}
In fact, the argument in Example~\ref{example:surface_bundles} is agnostic with
respect to dimension. Suppose $M$ is a closed hyperbolic manifold of dimension $n$
that fibers over the circle (by reason of Euler characteristic and Chern--Gauss--Bonnet
$n$ must be odd); cusped examples are known to exist in dimension $n=5$ and perhaps there are
examples in higher dimension 
(see \cite{Italiano_et_al_fiber}) and one might optimistically expect closed examples.
Reasoning exactly as above, any flow $X$ transverse to the foliation 
by fibers is quasigeodesic. Thus there are continuous maps $e^\pm:\H^n \to S^{n-1}_\infty$
taking each flowline of $\tilde{X}$ to its limit in forward and backward time. The images
$Z^\pm \subset S^{n-1}_\infty$ are disjoint, path-connected and $\pi_1(M)$-invariant;
however their higher connectivity (for instance, their homotopy groups) is a complete
mystery.
\end{example}

\begin{remark}
In the proof of Theorem~\ref{theorem:zippers_from_flows} the appeal to Frankel's theorem
is only used to show that certain quasigeodesic flows give rise to honest zippers
(as opposed to P-zippers); the argument that all quasigeodesic flows give rise to P-zippers
is entirely self-contained.
\end{remark}

\subsection{Pseudo-Anosov flows}

A flow $X$ on a 3-manifold $M$ is {\em Anosov} if it preserves a splitting of the tangent
bundle $TM=E^s\oplus E^u\oplus TX$ where $E^s$ resp. $E^u$ is uniformly expanded
resp. contracted by the flow. The bundles $E^s\oplus TX$ and $E^u \oplus TX$ are necessarily
integrable, and tangent to two 2-dimensional foliations. A flow is {\em pseudo-Anosov}
if it is Anosov away from finitely many singular orbits (possibly none), 
where the flow looks locally like the suspension of a $p$-prong singularity for some $p\ge 3$. 
Thus Anosov flows are special cases of pseudo-Anosov ones.

If $X$ is pseudo-Anosov, the leaf space of $\tilde{X}$ is a plane $P$,
and Fenley showed (\cite{Fenley_ideal_boundaries}, Thm.~A) that 
$P$ may be $\pi_1(M)$-equivariantly compactified to a closed disk by adding an ideal boundary circle 
(compare \cite{Calegari_Dunfield}, Thm.~3.8). 

Fenley (\cite{Fenley_ideal_boundaries}, \cite{Fenley_flow_circle}) has given necessary
and sufficient conditions for a pseudo-Anosov flow $X$ to be quasigeodesic. Pseudo-Anosov
flows without perfect fits are quasigeodesic, and more generally he shows that a pseudo-Anosov
flow is quasigeodesic if and only if it does not have chains of connected perfect fits
of unbounded length (\cite{Fenley_flow_circle}, Thm.~A). 
Thus, a pseudo-Anosov flow without perfect fits gives rise to a zipper, 
and a pseudo-Anosov flow with bounded chains of perfect fits gives rise to a P-zipper.

\subsection{Taut foliations}\label{subsection:taut_foliations}

Let $M$ be a hyperbolic 3-manifold, and let $\F$ be a taut co-orientable foliation of $M$.
Thurston sketched an argument \cite{Thurston_circles_II} 
that $\F$ should give rise to a universal circle $S^1_\un$ with a
$\pi_1(M)$-action, which bears a precise monotone relationship to the circles at infinity
of the leaves of $\tilde{\F}$. Although Thurston never provided complete
details of his construction, his proposal was carried out in \cite{Calegari_Dunfield} Thm.~6.2
(also see \cite{Calegari_foliations_book}, Chapter~7).
Actually there are at least {\em two} universal circles associated to $\F$ --- a leftmost
and a rightmost one --- and their relationship to the large scale geometry of
$\tilde{\F}$ is not at all straightforward.

Under some circumstances, one may construct a pseudo-Anosov flow $X$ transverse or
`almost transverse' to $\F$, and Fenley's criterion implies that these flows are quasigeodesic.
Specifically:
\begin{enumerate}
\item{when $\F$ is finite depth, Mosher \cite{Mosher_finite_depth_flow} constructed
a pseudo-Anosov flow $X$ almost transverse to $\F$ which Fenley--Mosher
\cite{Fenley_Mosher} later showed is quasigeodesic;}
\item{when $\F$ is $\R$-covered, Calegari \cite{Calegari_R_covered} and Fenley 
\cite{Fenley_R_covered} independently constructed a regulating pseudo-Anosov flow
$X$ transverse to $\F$ which Fenley \cite{Fenley_ideal_boundaries} later showed
has no perfect fits; and}
\item{when $\F$ has one-sided branching, Calegari \cite{Calegari_one_sided} constructed
a semi-regulating pseudo-Anosov flow $X$ transverse to $\F$ which
Fenley \cite{Fenley_ideal_boundaries} later showed has no perfect fits.}
\end{enumerate}
In the first case one obtains a P-zipper associated to $\F$, and in the other cases
one obtains a zipper. 

\begin{remark}
In fact, although most of the essential details are known to experts, Mosher's sequel to
\cite{Mosher_finite_depth_flow} in which the flow is really constructed is not
actually publically available. Michael Landry and Chi Cheuk Tsang are in the process of writing
up their own version of the construction; the first part is available as
\cite{Landry_Tsang}, and the second part is expected to be available soon.
Neither construction is straightforward.
For this and other reasons it would be valuable and interesting to give a direct
construction of zippers or P-zippers associated to a finite depth foliation.
\end{remark}

\begin{example}[Foliations with one-sided branching]
Suppose $M$ is hyperbolic and $\F$ is taut with one-sided branching. This means that
the leaf space $L$ of the universal cover $\tilde{\F}$ branches in the positive
direction (say) but not the negative direction. This means the following: if we define
a partial order $\prec$ on $L$ by $\mu \prec \lambda$ if there is a positively oriented
transversal to $\tilde{\F}$ from $\mu$ to $\lambda$, then for any two leaves $\lambda,\mu \in L$
there is $\nu \in L$ with $\nu \prec \lambda$ and $\nu \prec \mu$.

Let $\lambda$ be a leaf of $\tilde{\F}$.
In its intrinsic metric, $\lambda$ is quasi-isometric to the hyperbolic plane, and
therefore has a natural circle at infinity $S^1_\infty(\lambda)$.
Fenley showed \cite{Fenley_ideal_boundaries} that the natural inclusion 
$\lambda \to \tilde{M}=\H^3$ extends continuously to $S^1_\infty(\lambda) \to S^2_\infty$.
Let's denote the image by $b_\lambda$.

The following facts are true, and follow easily from the main results of 
\cite{Calegari_one_sided} and \cite{Fenley_ideal_boundaries}:
\begin{enumerate}
\item{for all $\lambda \in L$, the set $b_\lambda$ is compact, nonempty, path-connected,
and contains no embedded loop; and}
\item{if $\mu \prec \lambda$ then $b_\lambda \subset b_\mu$.}
\end{enumerate}
Thus, in its subspace topology, each $b_\lambda$ is a (compact) topological $\R$-tree.
No individual $b_\lambda$ is $\pi_1(M)$-invariant, but we may define
$b:=\cup_{\lambda \in L} b_\lambda$ and observe that $b$ is $\pi_1(M)$-invariant. 
We do not know if $b$ contains an embedded loop in its subspace topology; however,
we may give $b$ its weak topology as an increasing union, and in this topology it
is a (non-compact) topological $\R$-tree with a canonical planar embedding, 
and thereby gives rise to a circle with a $\pi_1(M)$-action.

What is the relation of $b$ to the zipper $Z^\pm$ associated to the quasigeodesic
pseudo-Anosov flow $X$ associated to $\F$?
\end{example}

\subsection{Flows from zippers}

Let $M$ be hyperbolic, and suppose that $Z^\pm$ is a minimal zipper for $M$ with universal
circle $S^1_\un$ and invariant laminations $\Lambda^\pm$.
Recent exciting work of Barthelm\'e--Bonatti--Mann \cite{Barthelme_Bonatti_Mann} and 
Frankel--Landry \cite{Frankel_Landry} gives conditions under which $\Lambda^\pm$ can
be completed to a pair of transverse singular foliations of the plane, and under
certain conditions (see \cite{Frankel_Landry}) 
one expects that this plane should be the leaf space of a
pseudo-Anosov quasigeodesic flow on $M$. Although there are many technical details
to work out in this program, it seems reasonable to conjecture a converse to the
main theorem of this section:
\begin{conjecture}
Let $M$ be a hyperbolic 3-manifold with a minimal P-zipper $Z^\pm$. Then there is a quasigeodesic
pseudo-Anosov flow $X$ on $M$ whose endpoint maps satisfy $e^\pm(P^\pm)=Z^\pm$. Furthermore, $X$
has no perfect fits if and only if $Z^\pm$ is a zipper.
\end{conjecture}

\section{Uniform quasimorphisms}\label{section:quasimorphisms}

A {\em quasimorphism} on a group is a real valued function that fails to be an additive
homomorphism in a controlled way (see Definition~\ref{definition:quasimorphism} 
for a precise definition). A quasimorphism is {\em uniform} if (roughly speaking) the
coarse level sets are coarsely connected.
The main theorem of this section (Theorem~\ref{theorem:zippers_from_quasimorphisms})
says that on any hyperbolic group, a uniform quasimorphism gives rise to a zipper.

\subsection{Quasimorphisms}

\begin{definition}[Quasimorphism]\label{definition:quasimorphism}
Let $G$ be a group. A function $\psi:G \to \R$ is a {\em quasimorphism} if there is some least
real constant $D(\psi)\ge 0$ (called the {\em defect}) so that for all $g,h\in G$ there is an
inequality $$|\psi(gh)-\psi(g)-\psi(h)|\le D(\psi)$$
\end{definition}

For an introduction to the theory of quasimorphisms and their properties, see e.g.\/
\cite{Calegari_scl}.

\begin{lemma}[Shrink and quantize]\label{lemma:shrink_and_quantize}
Let $\psi:G \to \R$ be a quasimorphism with defect $D(\psi)$. Let $t > 2D(\psi)$ and define
$\phi:G \to \Z$ by $\phi(g) := \lfloor \psi(g)/t + 1/2\rfloor$. Then $\phi$ is a quasimorphism with
$D(\phi)\le 1$.
\end{lemma}
\begin{proof}
For any real number $x$ the difference between $x$ and $\lfloor x+1/2\rfloor$ is at most $1/2$.
Since $D(\psi/t) < 1/2$ it follows from the triangle inequality that $D(\phi) < 2$. On the
other hand, the defect of $\phi$ is an integer, so $D(\phi) \le 1$.
\end{proof}

We say $\phi$ as above is obtained from $\psi$ by {\em shrink and quantize}.

\subsection{Uniform quasimorphisms}\label{subsection:uniform_quasimorphisms}

Let $G$ be a finitely generated group, with finite symmetric generating set $S$, and let $d$
be the word metric on $G$ induced by $S$. 

\begin{definition}
If $X\subset G$, a {\em $C$-chain} in $X$ is a sequence $x_0,x_1,\cdots,x_n$ of elements of $X$ with
$d(x_i,x_{i+1})\le C$ for all $i$.

A subset $X\subset G$ is {\em $C$-connected} for some integer $C>0$ if any two elements $x,y\in X$
may be joined by a $C$-chain in $X$.

Let $f:\N \to \N$ be some nondecreasing function. A subset $X\subset G$ is {\em $(C,f)$-connected} 
if it is $C$-connected, and if for all positive integers $m$ and for all $x,y\in X$ with $d(x,y)\le m$, 
there is a $C$-chain in $X$ from $x$ to $y$ of length at most $f(m)$.
\end{definition}

\begin{definition}\label{definition:uniform}
Let $C>0$ be an integer and $f:\N \to \N$ some nondecreasing function.
A quasimorphism $\phi:G \to \Z$ is {\em $(C,f)$-uniform} (or just {\em uniform} if $(C,f)$ are
left unspecified) if 
\begin{enumerate}
\item{the defect $D(\phi)\le 1$;}
\item{for any $g\in G$ there are $h^\pm$ with $d(g,h^\pm)\le C$ and $\phi(h^\pm)-\phi(g)=\pm 1$; and}
\item{for any integer $n$, the level set $\phi^{-1}(n)$ is $(C,f)$-connected.}
\end{enumerate}
\end{definition}

We refer informally to property 2 as the {\em wiggle property}.

\subsection{Equivalent definitions}

The following lemma shows that a quasimorphism with coarsely connected level sets may
be shrunk and quantized (in the sense of Lemma~\ref{lemma:shrink_and_quantize}) to
a quasimorphism that is uniform in the sense of Definition~\ref{definition:uniform}.

\begin{lemma}\label{lemma:coarsely_connected}
Let $\psi:G \to \R$ be a quasimorphism with defect $D(\psi)$. Suppose
\begin{enumerate}
\item{there is $g\in G$ for which $\psi(g)>2D(\psi)$;}
\item{there is $C_1>2\psi(g)$ and $C_2$ so that $\psi^{-1}[-C_1,C_1]$ is $C_2$-connected.}
\end{enumerate}
Let $t>2C_1$ and let $\phi:G \to \Z$ be obtained from $\psi$ by shrink and quantize
as in the proof of Lemma~\ref{lemma:shrink_and_quantize}. 
Then $\phi$ is $(C,f)$-uniform for some $(C,f)$.
\end{lemma}
\begin{proof}
Lemma~\ref{lemma:shrink_and_quantize} shows that $D(\phi)\le 1$, so that $\phi$ satisfies
the first property.

Let $g\in G$ be as in the statement of the lemma. For any positive integer $n$ we have 
$|\psi(g^n) - \psi(g^{n-1}) - \psi(g)| \le D(\psi)$. Since $\psi(g)>2D(\psi)$, by the
triangle inequality we have $\psi(g^n)> \psi(g^{n-1})+D(\psi)$ so by induction
$\psi(g^n)>nD(\psi)$ (and by a similar argument $\psi(g^{-n})<(1-n)D\psi$ for negative 
powers). On the other hand, the same inequality gives $\psi(g^n) < 2n\psi(g)$. It follows
that $\phi(g)=0$ and there is some constant $n_1$ so that for any $h$ there are $|n_2|,|n_3|\le n_1$
with $\phi(hg^{n_2}) = \phi(h)+1$ and $\phi(hg^{n_3}) = \phi(h)-1$. Thus $\phi$ satisfies
the wiggle property for some constant $C'$.

Finally we show the third property. Let $\phi(a)=\phi(b)=n$ and suppose $d(a,b)=m$. 
Then $|\phi(a^{-1}b)|\le 1$
and by the wiggle property we can find $c$ with $d(a^{-1}b,c)<C'$ and $\phi(c)=0$.
It follows that $|\psi(c)|<t/2$ so there is an integer $|k|<t/2D(\psi)$ so that 
if we write $c'=cg^k$ then $|\psi(c')|<C_1$. It follows from the hypothesis of the lemma
that $\id$ and $c'$ may be joined by a $C_2$-chain in $\psi^{-1}[-C_1,C_1] \subset \phi^{-1}(0)$. 
Now, $d(\id,c')\le d(a,b)+C'+|g|\cdot t/2D(\psi)$. Note that $C'+|g|\cdot t/2D(\psi)$ is a
constant independent of $d(a,b)$, so whatever $d(a,b)$ is there are only finitely many $c'$
within this distance of $\id$. Thus there is a uniform bound on the length of a
shortest $C_2$-chain $\id:=c_0,c_1,\cdots,c_n:=c'$ from $\id$ to $c'$ depending only on $d(a,b)=m$.
Define $f(m)$ to be this uniform bound plus 1. Then $a=ac_0,ac_1,\cdots,ac_n,b$ is a 
$C_2'$-chain from $a$ to $b$ of length $\le f(m)$ where $C_2':=\max(C_2,C'+|g|\cdot t/2D(\psi))$.
It is not necessarily true that $\phi(ac_i)=n$ for all $i$, but $|\phi(ac_i)-n|\le 1$ so by the
wiggle property we can find a new $C_2''$-chain $a=c_0',c_1',\cdots,c_n',b$ in $\phi^{-1}(n)$
of length at most $f(m)$ from $a$ to $b$, where $C_2'' = C_2'+2C'$.

Now define $C=\max(C',C_2'')$ and $f$ as above.
\end{proof}

\subsection{Staircases}

\begin{definition}
Let $\phi:G \to \Z$ be $(C,f)$-uniform. A {\em $(C,k)$-staircase} $a$ is an infinite 
sequence $a_0,a_1,a_2,\cdots$ of elements of $G$ for which
\begin{enumerate}
\item{$\phi(a_{i+1})=\phi(a_i)$ or $\phi(a_i)+1$ for all $i$;}
\item{$d(a_i,a_{i+1})\le C$ for all $i$; and}
\item{$\phi(a_{i+k})>\phi(a_i)$ for all $i$.}
\end{enumerate}
\end{definition}

\begin{lemma}\label{lemma:translate_is_staircase}
Let $\phi:G\to \Z$ be $(C,f)$-uniform. Then every $(C,k)$-staircase is $C'$-quasigeodesic for
$C'$ depending on $(C,f)$ and on $k$. Furthermore, if $G$ is hyperbolic then for every $k$ there is
a $k'$ so that if $a$ is a $(C,k)$-staircase and $g\in G$, there is a $(C,k')$-staircase
$a'$ with the same endpoint as $ga$.
\end{lemma}
\begin{proof}
A quasimorphism is Lipschitz in any word metric, so any $C$-chain on which
$\phi$ increases linearly is a quasigeodesic.

Define $a_i'=ga_i$. Since $\phi$ is a quasimorphism with $D(\phi)\le 1$
we have 
$$\phi(a_{i+3k}')\ge \phi(g) + \phi(a_{i+3k})-1 > \phi(g) + \phi(a_i)+1 \ge \phi(a_i')$$
Successively discard indices for which $\phi(a_{i+1}')<\phi(a_i')$ then
interpolate between successive $a_i'$ using properties 2 and 3 from 
Definition~\ref{definition:uniform}. This adjusts indices a uniformly bounded amount,
independent of $g$.
\end{proof}

\subsection{Blocks}

\begin{definition}[Block]
Let $\phi:G \to \Z$ be $(C,f)$-uniform. A {\em $(C,k)$-block} is a tuple $a,b,c_0,\cdots,c_n$ of
elements of $G$ for which
\begin{enumerate}
\item{$\phi(a)=\phi(b)$ and $\phi(c_i)=\phi(a)+1$ for all $i$;}
\item{$d(a,b)\le C$ and $d(c_i,c_{i+1})\le C$ for all $i$;}
\item{$d(a,c_0)\le C$ and $d(b,c_n)\le C$; and}
\item{$n\le k$.}
\end{enumerate}
\end{definition}
The pair $a,b$ is the {\em base} of the block and the $C$-chain $c_0,\dots,c_n$ is the
{\em roof} of the block.

\begin{lemma}[Extending block]\label{lemma:extending_block}
Let $a,b$ satisfy $\phi(a)=\phi(b)$ and $d(a,b)\le C$. Then there is a $k$ depending only
on $\phi$ so that $a,b$ is the base of a $(C,k)$-block. Furthermore if $c_0$ is any
element with $d(a,c_0)\le C$ and $\phi(c_0)=\phi(a)+1$ and/or $c_n$ is any element
with $d(b,c_n)\le C$ and $\phi(c_n)=\phi(b)+1$ then $a,b$ and/or $c_0$, $c_n$ may be
extended to a $(C,k)$-block with the same constant $k$ as above.
\end{lemma}
\begin{proof}
By the wiggle property we may obtain $c_0$ and/or $c_n$ with the desired properties and
then join $c_0$ to $c_n$ by a $C$-chain. The length of this chain depends only on the
distance from $c_0$ to $c_n$ which is at most $3C$.
\end{proof}

\begin{theorem}\label{theorem:zippers_from_quasimorphisms}
Let $G$ be a hyperbolic group and suppose that $G$ admits a uniform quasimorphism.
Then there is a zipper $Z^\pm \subset \partial_\infty G$ for $G$.
\end{theorem}
\begin{proof}
Let $\phi:G \to \Z$ be uniform.
By the definition of uniform, $\pm \phi:G \to \Z$ are $(C,f)$-uniform for some $C$ and
$f$. We shall show that there is a constant $k$, and $G$-invariant families
of $(C,k)$-staircases for $\pm \phi$
whose endpoints $Z^\pm \subset \partial_\infty G$ are a zipper.

By the wiggle property, some $(C,k)$-staircase exists.
By Lemma~\ref{lemma:translate_is_staircase}, if $a$ is a $(C,k)$-staircase 
then $ga$ is a $(C,k')$-staircase for some $k'$ independent of $g$. Thus once we
have one $(C,k)$-staircase, we can find a $G$-invariant family of them, at the cost of
adjusting the constant $k$. 
We shall show, for sufficiently large $k$, that the family $\PP^+$ of all $(C,k)$
staircases for $\phi$ has a path connected set of endpoints in $\partial_\infty G$.
Again by Lemma~\ref{lemma:translate_is_staircase} the set $G\PP^+$ of $G$-translates
of $\PP^+$ consists entirely of $(C,k)$ staircases for $\phi$ for some new $k$.
Furthermore its set of endpoints $Z^+$ is both $G$ invariant and a union of
path connected sets that pairwise intersect (since every $g\PP^+$ contains the entire
family $Ga$), and is therefore path connected.

Performing the same construction with $-\phi$ in place
of $\phi$ will produce another path connected and $G$-invariant subset 
$Z^- \subset \partial_\infty G$ disjoint from $Z^+$. This will prove the theorem.

Thus it suffices to show that any two $(C,k)$ staircases for $\phi$ can be included
in a family of $(C,k)$ staircases whose set of endpoints is path connected.
Let $a:=a_0,a_1,\dots$ and $b:=b_0,b_1,\dots$ be two $(C,k)$ staircases. 
By restricting to subsets and reordering we may assume $\phi(a_0)=\phi(b_0)=n$. 
Since $\phi$ is uniform, there is a $C$-chain in the level set $\phi^{-1}(n)$ from
$a_0$ to $b_0$, and by translation (say), every element of this $C$-chain is the 
initial element of some $(C,k)$-staircase; thus without loss of generality we may assume
$a$ and $b$ are adjacent elements of this chain, i.e.\/ $d(a_0,b_0)\le C$.

To prove the theorem we shall inductively apply Lemma~\ref{lemma:extending_block} to
interpolate between the staircases $a$ and $b$ by successively filling in blocks;
Figure~\ref{fill_in_blocks} indicates the argument.

\begin{figure}[htpb]
\labellist
\small\hair 2pt
\pinlabel $a_0$ at 385 15
\pinlabel $a_1$ at 370 110
\pinlabel $a_2$ at 340 110
\pinlabel $a_3$ at 310 110
\pinlabel $a_4$ at 250 205
\pinlabel $b_0$ at 515 15
\pinlabel $b_1$ at 530 110
\pinlabel $b_2$ at 570 110
\pinlabel $b_3$ at 615 205
\endlabellist
\centering
\includegraphics[scale=0.5]{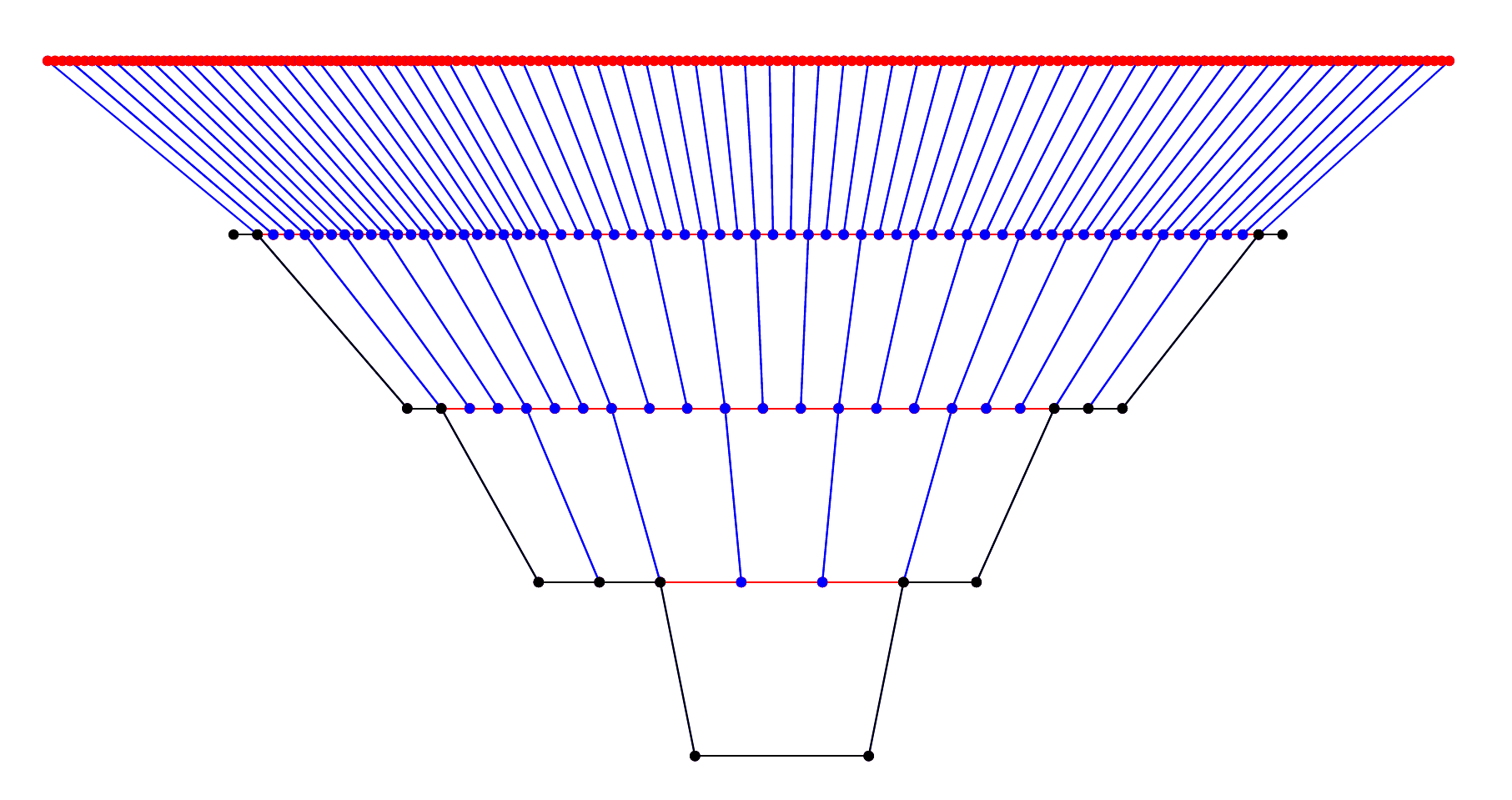}
\caption{The staircases $a$ and $b$ may be interpolated by successively filling in blocks.}
\label{fill_in_blocks}
\end{figure}

The staircases on the left and right of the figure (in black) are $a$ and $b$
respectively. Group elements at the same height have the same $\phi$ value, and every
segment has length $\le C$. The rest of the figure is filled in inductively as follows:
every horizontal segment of length $\le C$ may be extended vertically (in blue)
by the wiggle property, and then blocks may be filled in with roofs (in red)
by Lemma~\ref{lemma:extending_block}. The result of this inductive filling is a structure
as indicated in the figure. 

Consider the set of {\em all} $C$-chains in the figure that alternate between
vertical segments and $C$-chains of horizontal segments contained within a single roof.
Each such $C$-chain is a $(C,k)$ staircase by construction, and the set of endpoints
of this family in $\partial_\infty G$ is evidently path connected (one can think of
it as the image of the path that is the limit of the sequence of horizontal `paths'
joining an element of $a$ and an element of $b$ at each fixed height).
We have thus shown that the set of endpoints of the family $\PP^+$ of
all $(C,k)$ staircases for $\phi$ is path connected and the theorem is proved.
\end{proof}

\begin{remark}
The argument in Theorem~\ref{theorem:zippers_from_quasimorphisms} gives a practical
method to compute numerical approximations to zippers, assuming one has a relatively
concrete uniform quasimorphism $\phi$ to work with (e.g.\/ a homomorphism associated to
a fibration over $S^1$). One may construct an
approximate nontrivial path in $Z^+$ (say) by recursively extending blocks out to
some depth, and then transport this path around by the $G$-action to draw the rest of 
$Z^+$. This was how Figure~\ref{zipper_figure} was produced.
\end{remark}

\begin{example}[Slitherings]\label{example:slitherings}
A 3-manifold $M$ is said to {\em slither over $S^1$} if the universal cover $\tilde{M}$
fibers over $S^1$ in such a way that $\pi_1(M)$ acts by bundle automorphisms. The components
of the fibers are the leaves of a foliation $\tilde{\F}$ of $M$ covering a taut foliation
$\F$ of $M$. Slitherings are introduced in \cite{Thurston_circles_I} and many constructions
are given. For example, one may find examples where $M$ is a hyperbolic homology sphere.

Let $\rho:\tilde{M} \to S^1$ be the fibration, and $\tilde{\rho}:\tilde{M} \to \R$ a lift
to the universal cover. The point preimages of $\tilde{\rho}$ are exactly the leaves
of $\tilde{\F}$, so that we may identify the leaf space $L$ with $\R$. The action of
$\pi_1(M)$ on $\tilde{\F}$ induces therefore a (faithful orientation-preserving) action
of $\pi_1(M)$ on $\R$ by homeomorphisms, commuting with the group of integer translations
of $\R$. If we pick a basepoint $p \in \tilde{\rho}^{-1}(0) \in \tilde{M}$ the
function $\psi:\pi_1(M) \to \R$ defined by $\psi(g) =\tilde{\rho}(gp)$ is a quasimorphism.
Furthermore, the coarse level sets $\psi^{-1}[-C,C]$ are just the group elements $g$
for which $gp$ is between the leaves $\tilde{\rho}^{-1}(-C)$ and $\tilde{\rho}^{-1}(C)$;
when $C_1$ is big enough, this is evidently $C_2$-connected for some $C_2$, so that
$\psi$ satisfies the conditions of Lemma~\ref{lemma:coarsely_connected}. Thus $\phi$
obtained from $\psi$ by shrink and quantize is uniform, and gives rise to a zipper, and
consequently a universal circle.

Of course, $\F$ is $\R$-covered, so as pointed out in \S~\ref{subsection:taut_foliations} 
there is a regulating quasigeodesic flow $X$ which gives rise to a zipper as in
Theorem~\ref{theorem:zippers_from_flows}. Nevertheless, even in this case, the construction
of zippers and universal circles directly from the quasimorphism $\psi$ via
Theorem~\ref{theorem:zippers_from_quasimorphisms} is considerably more direct
than the arguments in \cite{Calegari_R_covered} or \cite{Fenley_R_covered}.
\end{example}

\subsection{Representations}

Quasimorphisms are very easy to construct, especially on hyperbolic manifolds.
If $M$ is a hyperbolic manifold and $\alpha$ is any 1-form at all, we may choose a
basepoint $p$ and define a map $\phi:\pi_1(M,p) \to \R$ by
$$\phi(g) = \int_{\gamma_g} \alpha$$
where $\gamma_g$ is the unique geodesic representative based at $p$ of the homotopy
class $g\in \pi_1(M,p)$. This is always a quasimorphism, with defect at most
$2\pi \|d\alpha\|_\infty$. However, there is no reason to expect in any particular
case that this quasimorphism is uniform.

Another class of quasimorphisms comes from group representations. Let
$\rho:\pi_1(M) \to G$ be a representation, where $G$ is either the group $\Homeo^+(S^1)$
or $\Sp(2n,\R)$ for some $n$. In either case $\pi_1(G)=\Z$ as a topological group, and
the obstruction to lifting $\rho$ to the universal covering group $\tilde{G}$ is
represented by an {\em Euler class} $e_\rho \in H^2(M;\Z)$ (in the case
$G=\Sp(2n,\R)$ this is usually called the {\em Maslov class}). Note that 
$G=\SL(2,\R)$ specializes both cases.

Suppose $e_\rho$ is zero, so that the representation lifts to 
$\tilde{\rho}:\pi_1(M) \to \tilde{G}$. In either case, $\tilde{G}$ admits a quasimorphism
$r:\tilde{G} \to \R$. In the case $G=\Homeo^+(S^1)$ this is the well-known (real-valued)
Poincar\'e rotation number; in the case $G=\Sp(2n,\R)$ this is the `symplectic
rotation number'. The composition $\phi = r\circ\tilde{\rho}$ is a quasimorphism on 
$\pi_1(M)$.

\begin{conjecture}
Let $M$ be a hyperbolic 3-manifold, let $\rho:G \to \SL(2,\C)$ be a lift of the
discrete faithful representation associated to the hyperbolic structure, and suppose
that $\rho$ has entries in $\SL(2,K)$ for some number field $K$ with a real Galois embedding
$\sigma:K \to K^\sigma \subset \R$. Let $\rho^\sigma:G \to \SL(2,K^\sigma) \subset \SL(2,\R)$ 
be the associated representation, and suppose $e_{\rho^\sigma} = 0$ in $H^2(M;\Z)$. Then
$\phi:=r\circ\tilde{\rho}^\sigma$ is uniform.
\end{conjecture}

\section{Uniform actions}\label{section:orders}

A {\em left order} on a group is a total order that is invariant under the (left) action
of the group on itself. A countable group is left orderable if and only if it acts faithfully
on $\R$ by orientation-preserving homeomorphisms. For such an action, the 
{\em up elements} of the group are the elements that move every point
in $\R$ in the positive direction. A faithful action is {\em uniform} if (roughly speaking) 
the partial order generated by up elements has upper bounds on a certain scale (and a similar
statement for the down elements). The main theorem of this section
(Theorem~\ref{theorem:zippers_from_orders}) says that on any hyperbolic group,
a uniform action gives rise to a zipper. 

\subsection{Up elements}

For the remainder of this section we fix a finitely generated left orderable 
group $G$ and a finite symmetric generating set $S$ for $G$.
For an introduction to left orderable groups and their basic properties, 
see e.g.\/ \cite{Navas_book}, especially Chapter~2.

As is well known, a (countable) group $G$ is left orderable if and only if 
there is a faithful action of $G$ on $\R$ by orientation-preserving homeomorphisms 
(\cite{Navas_book}, Theorem~2.2.19). If we conjugate this action by any
orientation-reversing homeomorphism of $\R$ we get the {\em opposite action},
related to an opposite reverse left order on $G$.
We would like to connect up the algebraic properties of the action to the geometry of $G$
and $M$; to do this we shall choose an action on $\R$ with good metric properties.

\begin{definition}
With notation as above, a orientation-preserving action of $G$ on $\R$ is {\em almost periodic} if
\begin{enumerate}
\item{there is a constant $K$ so that every $s\in S$ acts by a $K$-bilipschitz homeomorphism;}
\item{there is a constant $C$ so that for all $p\in \R$ and all $s\in S$ we have
$|s(p)-p|\le C$; and}
\item{for all $p\in \R$ there are elements $s^\pm\in S$ so that $s^+(p)-p\ge 1$ and
$s^-(p)-p \le -1$.}
\end{enumerate}
\end{definition}
Said in words, we can choose coordinates on $\R$ so that no generator or its inverse
stretches lengths too much, no generator moves any point too far, but every point is
moved a definite distance (up or down) by some generator.

The following is an amalgamation of \cite{Deroin} Theorem~3.1 and Lemma~4.2.
\begin{proposition}[Deroin]
Suppose $G$ is a finitely generated left-orderable group. Then $G$ admits a faithful
almost periodic action on $\R$.
\end{proposition}
Theorem~3.1 is the statement that any action may be conjugated to be bilipschitz, and
Lemma~4.2 says that any bilipschitz action may be further conjugated to satisfy the other
two properties. We remark for the benefit of the reader that the 
action produced by Deroin's proof of Theorem~3.1 has no global fixed points,
and therefore satisfies the hypotheses of Lemma~4.2.

Let's suppose we have a faithful almost periodic action of $G$ on $\R$. Define a
function $\phi:G \to \R$ by $\phi(g) = g^{-1}(0)$. From the definition of almost periodic
it follows that $\phi$ is Lipschitz.

\begin{definition}
An element $g\in G$ is {\em universally positive} (up) if there is a constant $C(g)$
so that $g^{-1}(p)-p\ge C(g)$ for all $p\in \R$.
\end{definition}

\begin{example}
For a minimal almost periodic action (in the sense of Deroin; see \cite{Deroin})
an element is an up element if and only if it acts freely and positively on $\R$.
Many naturally occurring left orders on 3-manifold groups (e.g.\/ typical orders
lifted from a circular order with trivial Euler class) have many up elements, and
it seems reasonable to study constructions that make use of them.

However such elements do not always exist!
As a concrete example, let $G=\langle a,b\rangle$ where $a:t \to 2t$ and $b$ has
compact support. For suitable choice of $b$ the action of $G$ has no global fixed
point, and we may always conjugate the action to be almost periodic. However,
for every element $g\in G$ either the exponent of $a$ (i.e.\/ the number of $a$s
in $g$ counted with sign) is zero in which case $g$ has compact support, or
the exponent is nonzero in which case $g$ or $g^{-1}$ takes a compact interval
properly inside itself. In either case $g$ has a fixed point.
\end{example}

Three manifold groups which are left orderable are often left orderable in many
different ways, and consequently admit many interesting actions on $\R$. Therefore we make the following:

\begin{conjecture}
Let $M$ be a closed 3-manifold. If $\pi_1(M)$ is left orderable, it admits a
faithful almost periodic action with up elements.
\end{conjecture}

\begin{example}
Any action on $\R$ with no global fixed points obtained by lifting an action 
on $S^1$ has many up elements. An $\R$-covered foliation is said to be
{\em uniform} if pairs of leaves in the universal cover are a bounded Hausdorff distance
apart; Thurston showed \cite{Thurston_circles_I} that such a foliation may always be
blown down to a slithering (compare with Example~\ref{example:slitherings}), and therefore the
action of the fundamental group on the leaf space has many up elements.
There are examples known \cite{Calegari_nonuniform}
of $\R$-covered foliations of hyperbolic 3-manifolds which are not uniform (in this sense),
but for these actions on $\R$ there are many up elements too.
\end{example}

\begin{definition}[Uniform order]\label{definition:uniform_order}
A faithful orientation-preserving action of a group $G$ on $\R$ is 
{\em up uniform} if it is almost periodic, and if there is some positive constant 
$C \ge 1$ so that for any two elements $x,y\in G$ with $d(x,y)\le 2C$ 
there is a $z\in G$ and two finite sequences of elements 
$x=x_0,\cdots,x_{m_1}=z$, $y=y_0,\cdots,y_{m_2}=z$ such that
\begin{enumerate}
\item{$d(x_i,x_{i+1})\le C$ and $d(y_i,y_{i+1})\le C$; and}
\item{each $x_i^{-1}x_{i+1}$ and $y_i^{-1}y_{i+1}$ is an up element.}
\end{enumerate}
An action is {\em uniform} if it is up uniform, and if the 
opposite action is also up uniform (possibly with different constants).
\end{definition}

Note that a group that admits a uniform action is left orderable.

\begin{remark}
Because there are only finitely many possibilities for $x^{-1}y$ with
$d(x,y)\le 2C$ one may impose a uniform bound on the sizes $m_1$, $m_2$ of the
up sequences in Definition~\ref{definition:uniform_order}. 
\end{remark}

\begin{remark}
It would be interesting to construct examples of almost periodic actions of groups
that are up uniform but not uniform. Perhaps this would be the analogue, in the world
of left orderable groups, of taut foliations with one-sided branching.
\end{remark}

\subsection{Zippers}

\begin{theorem}\label{theorem:zippers_from_orders}
Let $G$ be a hyperbolic group and suppose $G$ admits a uniform action on $\R$. 
Then there is a zipper $Z^\pm \subset \partial_\infty G$ for $G$.
\end{theorem}
\begin{proof}
Choose an almost periodic action on $\R$ which is up uniform for some constant $C$.
Let $U_C$ denote the set of up elements of
word length $\le C$. This set is nonempty, so for every
$g\in G$  we may construct an infinite sequence $g=a_0,a_1,\cdots$ where $a_i^{-1}a_{i+1}$
is in $U_C$; informally we call this an {\em up sequence}. 
Since $U_C$ is finite, by the definition of up elements
there is a positive constant $C_1>0$ so that $\phi(a_n) - \phi(g) \ge nC_1$ for
all $n$. Since $\phi$ is Lipschitz, this implies that every up sequence 
is uniformly quasigeodesic, independent of the choice of sequence. 

Let $Z^+\subset \partial_\infty G$ denote the set of endpoints of all 
up sequences starting at arbitrary $g\in G$, and let
$Z^-\subset \partial_\infty G$ be defined similarly with respect to up sequences for
the opposite action. Evidently both $Z^\pm$ are $G$-invariant and nonempty.
Furthermore if $a_0,a_1,\cdots$ is an up sequence for $Z^+$ and $b_0,b_1,\cdots$ is
an up sequence for $Z^-$ then $\phi(a_n) \to +\infty$ and $\phi(b_n)\to -\infty$;
since they are both quasigeodesic, and $\phi$ is Lipschitz, their endpoints are
disjoint. Thus $Z^+$ and $Z^-$ are disjoint.

It remains to show they are path-connected. We show this for $Z^+$.
Let $a_0,\cdots$ and $b_0,\cdots$ be up sequences with $d(a_0,b_0)\le C$. We
will show their endpoints are connected by a path in $Z^+$. Since there are up sequences
starting at any element of $G$ this is sufficient to prove the claim.

We shall construct an infinite directed planar graph $\Gamma$ in which every positively
oriented path is an up sequence. The directed graph will contain no directed loops,
so that it defines a partial order on its vertices, which are elements of $G$.
The `left' and `right' edges of $\Gamma$ are
the sequences $a_0,\cdots$ and $b_0,\cdots$. By the definition of up uniform, there are
finite up sequences $a_0=x_0,x_1,\cdots,x_{m_1}$ and $b_0=y_0,y_1,\cdots,y_{m_2}$ with
$x_{m_1}=y_{m_2}$; attaching these directed segments makes $\Gamma$ into a `W' shape;
we refer to this operation as {\em attaching a $\wedge$}. 

Let us call a $\vee$ in $\Gamma$ a triple of elements $a,b,c$ where $b$ is joined to
both $a$ and $c$ by a positively oriented edge, but $a$ and $c$ have no common
upper bound; this means that $b^{-1}a$ and
$b^{-1}c$ are both up elements of length $\le C$. In particular, $d(a,c)\le 2C$
so we may attach a $\wedge$ to the top of this $\vee$. Do this inductively, creating
a structure as in Figure~\ref{upper_graph}.

\begin{figure}[htpb]
\centering
\includegraphics[scale=0.3]{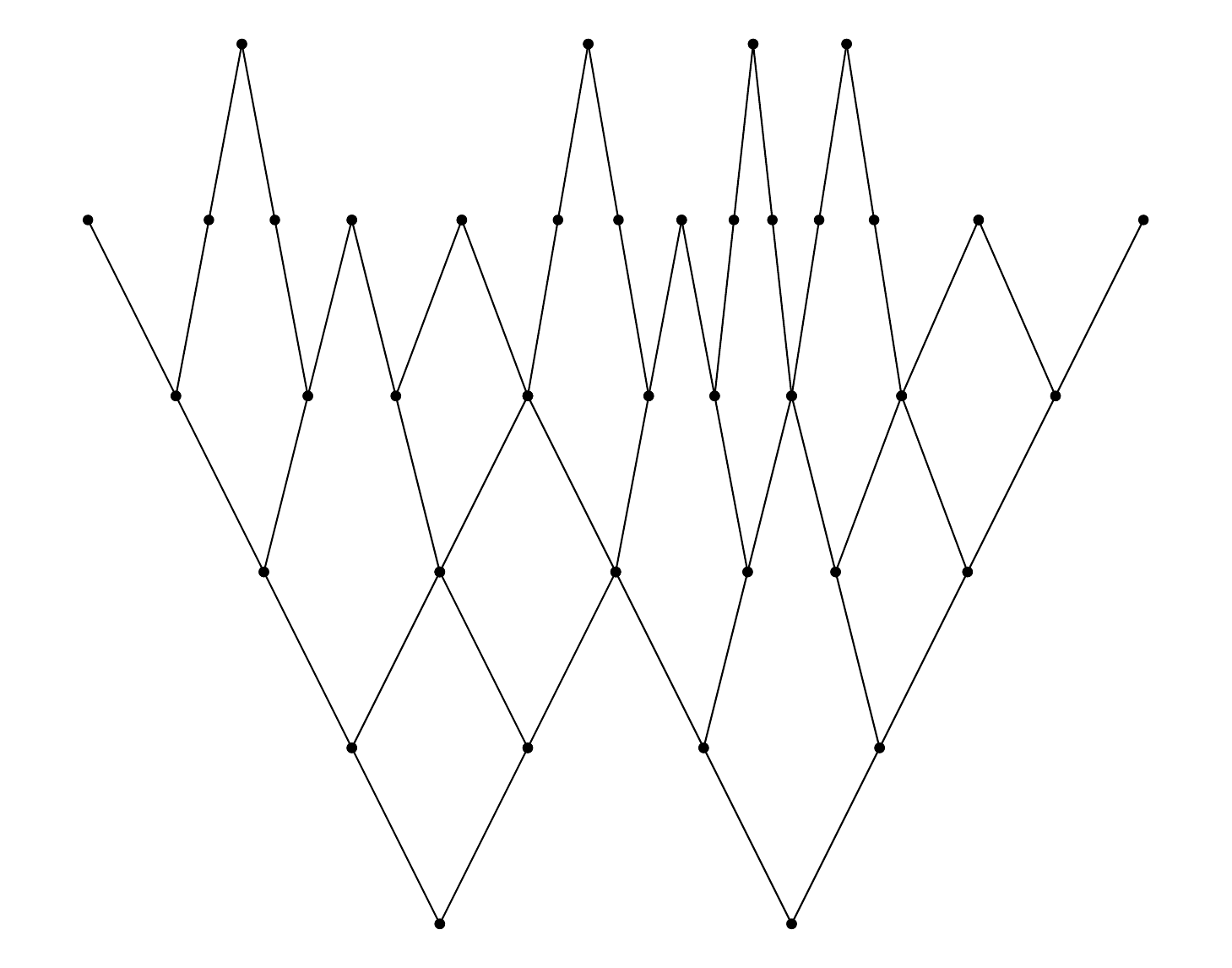}
\caption{$\Gamma$ is built inductively from $W$ by successively capping every 
$\vee$ with a $\wedge$.}
\label{upper_graph}
\end{figure}

The obvious planar structure defines a partial order on the set of infinite directed
ascending paths, where $\gamma_1 \ge \gamma_2$ if $\gamma_1$ is nowhere to the left of 
$\gamma_2$. A maximal ordered subset of such paths is a Cantor set, and the gaps in this
Cantor set are pairs of paths with infinitely many vertices in common, which have
a common endpoint at infinity. Thus the (continuous!) map from paths to endpoints
factors through a surjective map from a Cantor set to an interval, and the
set of directed endpoints of $\Gamma$ is a path-connected subset of $\partial_\infty G$.
\end{proof}

\subsection{Orders and quasimorphisms}

It is tempting to try to connect the concepts of uniform actions and uniform quasimorphisms.
Here is one potential point of contact. 

Let $G$ be a hyperbolic group, and let's fix a faithful orientation-preserving
action of $G$ on $\R$. If $S \subset G$ is any (not necessarily symmetric!) finite
subset that generates $G$ as a semigroup 
we may define a quasimorphism $\phi_S$ by {\em antisymmetrizing} the $S$ word metric.
That is, $\phi_S(g) = |g|_S - |g^{-1}|_S$. 
Let's fix a symmetric generating set $S$ for $G$, and for any $n$ define $S_n$ 
to be the union of $S$ with all elements $g\in G$ of word length $\le n$ for which $g(0)\le 0$,
and then construct the associated quasimorphism $\phi_n: = \phi_{S_n}$.

\begin{question}
What conditions on the action ensure that $\phi_n$ is a uniform quasimorphism for some $n$? 
If the given action is uniform, is $\phi_n$ necessarily uniform for some $n$?
\end{question}

\section{Acknowledgements}

We would like to thank Nathan Dunfield, Sergio Fenley,
Steven Frankel, KyeongRo Kim, Michael Landry, Curt McMullen and the anonymous referee
for valuable comments, assistance and encouragement.

\end{document}